\documentclass[12pt]{amsart}
  \usepackage{latexsym}
  \usepackage[all,2cell]{xy}
  \usepackage{amsfonts}
  \usepackage{amsthm}
  \usepackage{amsmath}
  \usepackage{amssymb}
 \numberwithin{equation}{section}
 \usepackage{mathptmx}
  \usepackage{hyperref}
\usepackage{multicol}
\usepackage{array}

\UseAllTwocells
\CompileMatrices

\def\ox{\otimes}
\def\x{\times}
\def\bx{\boxtimes}
\def\1c#1{\stackrel{#1}{\to}}
\def\Mnd{\mathsf{Mnd}\,}
\def\C{\mathcal C}
\def\M{\ensuremath{\mathcal M}}
\def\N{\mathcal N}
\def\Cat{\mathsf{Cat}}
\def\wbm{\mathsf{Wbm}}
\def\fsbm{\mathsf{Sfbm}}
\def\can{\mathsf{can}}
\def\runit{\varrho}
\def\lunit{\lambda}
\def\mox{\raise.1cm\hbox{${}_\ox$}}

\newcommand{\ot}{\otimes}
\renewcommand{\t}{\times}

\renewcommand{\epsilon}{\varepsilon}
\renewcommand{\phi}{\varphi}

  \newtheorem{proposition}{Proposition}[section]
  \newtheorem{lemma}[proposition]{Lemma}

  \newtheorem{theorem}[proposition]{Theorem}

  \theoremstyle{definition}
  \newtheorem{definition}[proposition]{Definition}
  \newtheorem{example}[proposition]{Example}

  \theoremstyle{remark}
  \newtheorem{remark}[proposition]{Remark}

  \newcounter{c}
  
  \newcommand{\etyk}[1]{\vspace{-7.4mm}$$\begin{equation}\Label{#1}
  \addtocounter{c}{1}}
  \renewcommand{\]}{\ifnum \value{c}=1 $$\else \end{equation}\fi}
\begin{document}

 \title{Weak bimonads and weak Hopf monads}

 \author{Gabriella B\"ohm}
 \address{Research Institute for Particle and Nuclear Physics, Budapest,
 \newline\indent 
H-1525 Budapest 114, P.O.B.\ 49, Hungary}
\email{G.Bohm@rmki.kfki.hu}

\author{Stephen Lack}
\address{School of Computing and Mathematics University of Western 
Sydney,
 \newline\indent 
Locked Bag 1797 Penrith South DC NSW 1797, Australia and 
\newline \indent 
Mathematics Department Macquarie University, NSW 2109 Australia.}
\email{s.lack@uws.edu.au; steve.lack@mq.edu.au}

\author{Ross Street}
\address{Mathematics Department Macquarie University, NSW 2109 Australia.}
\email{ross.street@mq.edu.au}

\date{February 2010}

\begin{abstract}
We define a {\em weak bimonad} as a monad $T$ on a monoidal category
${\mathcal M}$ with the property that the Eilenberg-Moore category $\M^T$ is 
monoidal and the forgetful functor $\M^T \to \M$ is
separable Frobenius. 
Whenever $\M$ is also Cauchy complete, a simple set of axioms is provided,
that characterizes the monoidal structure of $\M^T$ as a weak lifting of the
monoidal structure of $\M$.
The relation to bimonads, and the relation to weak bimonoids in a braided
monoidal category are revealed. 
We also discuss antipodes, obtaining the notion of weak Hopf monad.
\end{abstract}
  
\maketitle


\section*{Introduction}

{\em Bialgebras} (say, over a field) have several equivalent
characterizations. One of the most elegant ones is due to Pareigis, who proved
that an algebra $A$ over a field $K$ is a bialgebra if and only if the
category of (left or right) $A$-modules is monoidal and the forgetful functor
from the category of $A$-modules to the category of $K$-vector spaces is
strict monoidal. This fact extends to bialgebras in any braided monoidal
category \cite{Maj:br_Hopf}.  

Pareigis' characterization of a bialgebra was the starting point of Moerdijk's
generalization in \cite{Moe:Hopf_monad} of bialgebras to monoidal categories
possibly without a braiding. He defined a {\em bimonad} (originally called a
Hopf monad) as a monad $T$ on a monoidal category $\M$, such that
the Eilenberg-Moore category $\M^T$ of $T$-algebras is monoidal and
the forgetful functor $\M^T \to \M$ is strict monoidal. That is, the monoidal
structure of $\M$ lifts to $\M^T$. Because liftings of functors (respectively,
of natural transformations) are described by 1-cells (respectively, by
2-cells) in the 2-category $\Mnd(\Cat)$ of monads (in the notation of
\cite{Str:FTM}), Moerdijk's definition says that a monad is a bimonad if
and only if the functor induced by the monoidal unit of ${\mathcal M}$, from
the terminal category to $\M$, and the functor provided by the monoidal
product of $\M$, from $\M \x {\mathcal M}$ to $\M$, both admit the structure
of a 1-cell in $\Mnd(\Cat)$, and the coherence natural isomorphisms in $\M$
are 2-cells in $\Mnd(\Cat)$.  
In \cite{McCr:opmon_mon}, McCrudden showed that a bimonad is the same as an
opmonoidal monad, that is, a monad in the 2-category of monoidal categories,
opmonoidal functors and opmonoidal natural transformations. Equivalently,
bimonads are the same as monoids in a multicategory of monads on a monoidal
category. 

Pareigis' characterization of a bialgebra was generalized to a {\em weak
bialgebra} \cite{Nill:WBA}, \cite{BNSz:WHA_I} by Szlach\'anyi in 
\cite{Szl:Bru}. He proved that an algebra $A$ over a field $K$ is a weak
bialgebra if and only if the category of (left or right) $A$-modules is
monoidal and the forgetful functor from the category of $A$-modules to the
category of $K$-vector spaces obeys the so called {\em separable Frobenius}
condition. The latter means that the forgetful functor admits both a monoidal
and an opmonoidal structure that satisfy some compatibility relations: see
Definition \ref{def:Frob-sep}. 
These (op)monoidal structures are no longer strict. In particular, the monoidal
unit of the category of $A$-modules is not $K$ as a vector space but a
non-trivial retract of $A$. Also, the monoidal product of two $A$-modules is
not their $K$-module tensor product but a linear retract of it. 

Weak bialgebras can be defined in any braided monoidal category, see
\cite{AAatal:weak_cleft} and \cite{PaSt:WH_mon}, as objects possessing both a
monoid and a comonoid structure, subject to compatibility axioms that
generalize those in \cite{Nill:WBA} and \cite{BNSz:WHA_I} in the case of a
symmetric monoidal category of vector spaces. The resulting category of
modules was investigated in \cite{PaSt:WH_mon}.

The aim of this paper is to generalize weak bialgebras to monoidal categories
possibly without a braiding. Inspired by Szlach\'anyi's characterization of a
weak bialgebra, we define a {\em weak bimonad} as a monad $T$ on a monoidal
category $\M$, with extra structure making $\M^T$ monoidal and the forgetful
functor $\M^T\to\M$ separable Frobenius.

For a weak bimonad $T$, the forgetful functor $\M^T \to \M$ is no
longer strict monoidal, hence the monoidal structure of the domain category
$\M$ does not lift to $\M^T$, and so the monoidal unit and the monoidal
product of $\M$ are no longer 1-cells in $\Mnd(\Cat)$. 
However, the notion of lifting ${\overline F}: \M^T \to {\mathcal
M}^{\prime T'}$ of a functor $F:\M \to \M'$ was weakened
in \cite{Bohm:weak_th_mon} by replacing commutativity of the diagram of
functors   
$$
\xymatrix{
\M^T \ar[r]^-{\overline F} \ar[d]_-{U}& 
\M^{\prime T'}\ar[d]^-{U'}\\
\M \ar[r]^-{F} &
\M'
}
$$
by the existence of a split natural monomorphism $i:U' {\overline F} \to F
U$. A weak lifting of a natural transformation 
is defined as a natural transformation between the lifted functors that
commutes with the natural monomorphisms $i$ in the evident sense.  
Weak liftings of functors and of natural transformations in a locally Cauchy
complete 2-subcategory of $\Cat$, are related to 1-cells and 2-cells
in a 2-category $\Mnd^i(\Cat)$ in \cite{Bohm:weak_th_mon}, extending
$\Mnd(\Cat)$. 

In Section \ref{sec:ax} we give an interpretation of the axioms of a weak
bimonad (on a Cauchy complete monoidal category), 
similar to the interpretation of a bimonad in \cite{Moe:Hopf_monad}.  
While for a bimonad $T$ the monoidal structure of $\M^T$ 
is given by lifting of the monoidal structure in the domain category $\M$, for
a weak bimonad $T$ the monoidal product in $\M^T$ is a weak lifting of the
monoidal product in $\M$, the monoidal unit is a weak 
lifting of the functor $ 1\to \M \stackrel{T}{\to} \M$,
the associativity
constraint is a weak lifting of the associativity constraint in $\M$ and the
unit constraints are weak liftings of certain morphisms in $\M$ constructed
from the other data.  

By results in \cite{PaSt:WH_mon}, a weak bimonoid in a braided monoidal
category can be described as a quantum category over a separable Frobenius
base monoid. Extending this result in Section \ref{sec:Frob_sep_base}, we 
establish an equivalence between the category of weak bimonads on a
Cauchy complete monoidal category $\M$ and the category of bimonads on
bimodule categories over separable Frobenius monoids in $\M$. 

In Section \ref{sec:weak_bimonoid} we show that weak bimonoids in a braided
monoidal category (cf. \cite{PaSt:WH_mon}, \cite{AAatal:weak_cleft}) induce
weak 
bimonads. In certain braided monoidal categories the converse can also be
proved: if a monoid induces a weak bimonad then it admits the structure of a
weak bimonoid. 

In Section \ref{sec:antipode}, using the result in Section
\ref{sec:Frob_sep_base} that any weak bimonad (on a Cauchy complete monoidal
category) can be regarded as a bimonad (on another monoidal category), we
define a {\em weak right Hopf monad} to be a weak bimonad such that the
associated bimonad is a right Hopf monad in the sense of \cite{BLV} and
\cite{ChiklackStr:Homf_mon}; there is a companion result involving weak left
Hopf monads and left Hopf monads. 
A weak bimonoid in a Cauchy complete braided monoidal category is shown to
induce a weak right Hopf monad by tensoring with it on the right if and only
if it is a weak Hopf monoid in the sense of \cite{AAatal:weak_cleft} and
\cite{PaSt:WH_mon}; once again there is a companion result with left in place
of right. 

\subsection*{Notation and conventions.} The monoidal categories in this paper
are not necessarily strict but in order to simplify our expressions, we omit 
explicit mention of their coherence isomorphisms wherever possible. 

Recall that, for an opmonoidal functor $F:(\N,\bx,R)\to(\M,\ot,K)$
with opmonoidal structure $i_{X,Y}:F(X\bx Y)\to FX\ot FY$ and 
$i_0:FR\to K$, the diagram 
$$\xymatrix{
F(X\bx Y\bx Z) \ar[r]^{i_{X\bx Y,Z}} \ar[d]_{i_{X,Y\bx Z}} & 
F(X\bx Y)\bx FZ \ar[d]^{i_{X,Y}\bx FZ} \\
FX\bx F(Y\bx Z) \ar[r]_{FX\bx i_{Y,Z}} & FX\bx FY\bx FZ }$$
commutes. We sometimes write $i^{(3)}_{X,Y,Z}$ for the common composite,
and we use an analogous notation for monoidal functors.

We say that a category is {\em Cauchy complete} provided that idempotent 
morphisms in it split. 


\section{Weak bimonads and their Eilenberg-Moore category of algebras}
\label{sec:ax}

The definition of weak bimonad is based on the notion of 
separable Frobenius functor introduced in \cite{Szl:Bru}:

\begin{definition}\label{def:Frob-sep}
A functor $F$ from a monoidal category $({\mathcal N},\bx,R)$ to a
monoidal category $(\M,\ox,K)$ is said to be {\em
separable Frobenius} when it is equipped with a monoidal structure
$p_{X,Y}: FX \ox FY \to F(X\bx Y)$, $p_0:K \to FR$ and an
opmonoidal structure $i_{X,Y}: F(X\bx Y) \to FX \ox FY$,
$i_0: FR \to K$ such that, for all objects $X,Y,Z$ in ${\mathcal N}\!\!\!$,
the following diagrams commute:
$$\xymatrix{
FX\ot F(Y\bx Z) \ar[r]^{FX\ot i_{Y,Z}} \ar[d]_{p_{X,Y\bx Z}} & 
FX\ot FY\ot FZ \ar[d]^{p_{X,Y}\ot FZ} \\
F(X\bx Y\bx Z) \ar[r]_{i_{X\bx Y,Z}} & F(X\bx Y)\ot FZ
}\qquad 
\xymatrix{
F(X\bx Y)\ot FZ \ar[r]^{i_{X,Y}\ot FZ} \ar[d]_{p_{X\bx Y,Z}} & 
FX\ot FY\ot FZ \ar[d]^{FX\ot p_{Y,Z}} \\
F(X\bx Y\bx Z) \ar[r]_{i_{X,Y\bx Z}} & FX\ot F(Y\bx Z) }
$$
$$\xymatrix @R1pc {
& FX\ot FY \ar[dr]^{p_{X,Y}} \\
F(X\bx Y) \ar[ur]^{i_{X,Y}} \ar@{=}[rr] && F(X\bx Y) }$$
\end{definition}

\begin{example}\label{ex:Frob-sep}
(1) Strong monoidal functors are clearly separable Frobenius. 

(2) The composite of separable Frobenius functors is separable Frobenius,
cf. \cite{DaPa:Frob-sep}.

(3) In a monoidal category $(\M,\ox,K)$ possessing (appropriate)
coequalizers preserved by $\ox$, one may consider the monoidal
category ${}_R \M_R$ of bimodules over a monoid $R$ in $\M$. The monoidal 
product is provided by the $R$-module tensor product and the monoidal unit is
$R$. Justifying the terminology, the 
forgetful functor ${}_R \M_R\to \M$ is separable Frobenius if and only if $R$
is a separable Frobenius monoid; that is, a Frobenius monoid in the sense of
\cite{Str:Frob_mon} such that, in addition, composing its comultiplication
$R\to R \ox R$ with its multiplication $R\ox R \to R$ yields the identity
morphism $R$. 
\end{example}

\begin{definition}\label{def:wbm}
A {\em weak bimonad} on a monoidal category $(\M,\ox,K)$ is a monad
$(T,m,u)$ on $\M$ equipped with a monoidal structure on the Eilenberg-Moore
category ${\mathcal M}^T$ and a separable Frobenius structure on the forgetful
functor ${\mathcal M}^T\to \M$.  
\end{definition}

The main aim of this section is to find an equivalent formulation of
Definition \ref{def:wbm} -- in the spirit of the descriptions of bimonads in
\cite{Moe:Hopf_monad} and \cite{McCr:opmon_mon} (there called Hopf monads).
\medskip

If a monad $T$ possesses a monoidal Eilenberg-Moore category $({\mathcal
  M}^T,\bx, (R,r))$ then, for any $T$-algebras $(A,a)$ and $(B,b)$,
there is a $T$-algebra $(A,a)\bx (B,b)$ that we denote by $(A \Box B,
a\Box b)$. (Note that by definition $a \Box b$ is a morphism $T(A\Box B) \to A
\Box B$ in  
$\M$, while \hbox{$a \bx b$} is a morphism $(TA,m_A)\bx (TB,m_B)
\to (A,a) \bx (B,b)$ in $\M^T$; that is, a morphism $TA \Box TB
\to A \Box B$ in $\M$.  Note also that $A\Box B$ depends not just on $A$ and
$B$ but on the algebras $(A,a)$ and $(B,b)$.)

In order to get started, we need the following basic observation:

\begin{proposition}\label{prop:T_op-mon}
Consider a monad $(T,m,u)$ on a monoidal category $({\mathcal
  M},\ox,K)$ equipped with a monoidal Eilenberg-Moore category $({\mathcal
  M}^T,\bx, (R,r))$. If the forgetful functor $U:\M^T \to \M$ admits both a
monoidal structure $(p,p_0)$ and an opmonoidal structure $(i,i_0)$ then $T$
is opmonoidal, with  $\tau_0$ and $\tau_{X,Y}$ given, respectively,  by the
composite morphisms  
\begin{equation}\label{eq:tau_0}
\xymatrix{
TK \ar[r]^-{Tp_0}&
TR \ar[r]^-{r}&
R \ar[r]^-{i_0}&
K
}
\qquad \textrm{and}
\end{equation}
\begin{equation}\label{eq:tau_XY}
\xymatrix{
T(X \ox Y) \ar[r]^-{T(u_X \ox u_Y)}&
T(TX \ox TY) \ar[r]^-{Tp_{TX,TY}}&
T(TX \Box TY) \ar[r]^-{m_X\Box m_Y}&
TX \Box TY \ar[r]^-{i_{TX,TY}}&
TX \ox TY.
}
\end{equation}
\end{proposition}

\begin{proof}
Since $U$ is monoidal, its left adjoint $F$ is opmonoidal. Since $U$ is also
opmonoidal, so is $T=UF$. The explicit form of the structure morphisms
\eqref{eq:tau_0} and \eqref{eq:tau_XY} is immediate. 
\end{proof}

At this point we can now state one characterization of weak bimonads:

\begin{theorem}\label{thm:main}
Let $T=(T,m,u)$ be a monad on a monoidal category $(\M,\ot,K)$ in which
idempotents split. To give  
$T$ the structure of a weak bimonad is equivalently to give the endofunctor
$T$ the structure of an opmonoidal functor $(T,\tau,\tau_0)$ in such a way 
that the following conditions hold:
\begin{equation}
  \label{eq:tau_1}
  \xymatrix{
T^2(X\ot TK) \ar[r]^-{T\tau_{X,TK}} & T(TX\ot T^2K) \ar[r]^{T(TX\ot m_K)} & 
T(TX\ot TK) \ar[r]^-{T(TX\ot\tau_0)} & T^2 X \ar[d]^{m_X} \\
T(X\ot TK) \ar[u]^{Tu_{X\ot TK}} \ar[r]_{\tau_{X,TK}} & TX\ot T^2 K
\ar[r]_{TX\ot m_K} & 
TX\ot TK \ar[r]_-{TX\ot\tau_0} & TX}
\end{equation}
\begin{equation}
  \label{eq:tau_2}
  \xymatrix{
T^2(TK\ot X) \ar[r]^{T\tau_{TK,X}} & T(T^2 K\ot TX) \ar[r]^{T(m_K\ot TX)} & 
T(TK\ot TX) \ar[r]^-{T(\tau_0\ot TX)} & T^2 X \ar[d]^{m_X} \\
T(TK\ot X) \ar[u]^{Tu_{TK\ot X}} \ar[r]_{\tau_{TK,X}} & T^2 K\ot TX
\ar[r]_{m_K\ot TX} & 
TK\ot TX \ar[r]_-{\tau_0\ot TX} & TX}
\end{equation}
\begin{equation}
  \label{eq:tau_3}
  \xymatrix @C 1pc {
X\ot T(Y\ot Z) \ar[rr]^{X\ot \tau_{Y,Z}} && X\ot TY\ot TZ \ar[rrr]^{u_{X\ot TY}\ot TZ} &&&
T(X\ot TY)\ot TZ \ar[rr]^-{\tau_{X,TY}\ot TZ} && TX\ot T^2 Y\ot TZ \ar[d]^{TX\ot
  m_Y\ot TZ} \\ 
X\ot Y\ot Z \ar[u]^{X\ot u_{Y\ot Z}} \ar[rr]_{u_{X\ot Y\ot Z}} && T(X\ot Y\ot Z)
\ar[rrr]_{\tau_{X\ot Y,Z}} &&&
T(X\ot Y)\ot TZ \ar[rr]_-{\tau_{X,Y}\ot TZ} && TX\ot TY\ot TZ }
\end{equation}
\begin{equation}
  \label{eq:tau_4}
  \xymatrix @C 1pc{
T(X\ot Y)\ot Z \ar[rr]^{\tau_{X,Y}\ot Z} && TX\ot TY\ot Z \ar[rrr]^{TX\ot u_{TY\ot Z}} &&&
TX\ot T(TY\ot Z) \ar[rr]^-{TX\ot \tau_{TY,Z}} && TX\ot T^2 Y\ot TZ \ar[d]^{TX\ot
  m_Y\ot TZ} \\ 
X\ot Y\ot Z \ar[u]^{u_{X\ot Y}\ot Z} \ar[rr]_{u_{X\ot Y\ot Z}} && T(X\ot Y\ot Z)
\ar[rrr]_{\tau_{X\ot Y,Z}} &&& 
T(X\ot Y)\ot TZ \ar[rr]_{\tau_{X,Y}\ot TZ} && TX\ot TY\ot TZ }
\end{equation}
\begin{equation}
  \label{eq:tau_5}
  \xymatrix{
T^2(X\ot Y) \ar[r]^{T\tau_{X,Y}} \ar[d]_{m_{X\ot Y}} & T(TX\ot TY)
\ar[r]^{\tau_{TX,TY}} & 
T^2X\ot T^2 Y \ar[d]^{m_X\ot m_Y} \\
T(X\ot Y) \ar[rr]_{\tau_{X,Y}} && TX\ot TY }
\end{equation}
\end{theorem}

We shall spend the rest of the section proving this theorem as well as
formulating a further characterization in terms of weak lifting.
One half of the theorem we prove immediately:

\begin{proposition}
For any weak bimonad, equations~\eqref{eq:tau_1}--\eqref{eq:tau_5} hold when
the endofunctor is given the opmonoidal structure of 
Proposition~\ref{prop:T_op-mon}.
\end{proposition}

\proof
$(R,r)$ is the monoidal unit in $\M^T$ and the coherence natural isomorphisms
in $\M^T$ are $T$-algebra morphisms.
For any morphisms $f:(A,a)\to  (A',a')$ and $g:(B,b)\to (B',b')$ of
$T$-algebras, $f\bx g$ is a morphism of $T$-algebras, so that
\begin{equation}
\label{eq:c_3}
\xymatrix{
T(A\Box B) \ar[d]_{a\Box b} \ar[r]^{T(f\bx g)} & T(A'\Box B') \ar[d]^{a'\Box b'} \\
A\Box B \ar[r]_{f\bx g} & A'\Box B' }
\end{equation}
commutes. 
By \eqref{eq:tau_XY} and unitality of the $T$-action $mX \Box mY$, 
the diagram 
\begin{equation}
\label{eq:tau&u}
\xymatrix{
X\ot Y \ar[d]_{u_{X\ot Y}}  \ar[r]^{u_X\ot u_Y} & TX\ot TY \ar[r]^{p_{TX,TY}} & 
TX\Box TY \ar[d]^{i_{TX,TY}} \\
T(X\ot Y) \ar[rr]_{\tau_{X,Y}} && TX\ot TY
}
\end{equation}
commutes, for any objects $X,Y$ of $\M$.
Hence a straightforward computation, using these facts together with the
unitality of $m$ and with the opmonoidality of $(U,i,i_0)$,
shows that both routes around \eqref{eq:tau_1} are equal to 
$$\xymatrix{
T(X\ot TK) \ar[r]^{T(X\ot Tp_0)} & T(X\ot TR) \ar[r]^{T(X\ot r)} & 
T(X\ot R) \ar[r]^{T(u_X\ot R)} & T(TX\ot R) \ar[r]^-{Tp_{TX,R}} & 
T^2 X \ar[r]^{m_X} & TX. }$$
Equality
\eqref{eq:tau_2} is proved symmetrically. 
In view of \eqref{eq:tau&u}, the bottom path of \eqref{eq:tau_3} is equal to 
$$\xymatrix @C 3pc {
X\ot Y\ot Z \ar[r]^-{u_X\ot u_Y\ot u_Z} &  TX\ot TY\ot TZ
\ar[r]^{p^{(3)}_{TX,TY,TZ}} & 
TX\Box TY\Box TZ \ar[r]^{i^{(3)}_{TX,TY,TZ}} & TX\ot TY\ot TZ .}$$
This expression is checked to be equal also to the upper path of
\eqref{eq:tau_3}, by applying \eqref{eq:tau&u} repeatedly, and using
monoidality of $(U,p,p_0)$ and the first property in
Definition \ref{def:Frob-sep} of the separable Frobenius functor $U$. Equality 
\eqref{eq:tau_4} is proved symmetrically, using the second property in
Definition \ref{def:Frob-sep} of the separable Frobenius functor $U$ instead
of the first one. 
Finally, by \eqref{eq:tau_XY}, by naturality of $i$ and
$p$, by \eqref{eq:c_3}, and by unitality of $m$, we deduce that
$(m_X \ox m_Y)\circ \tau_{TX,TY}= i_{TX,TY}\circ (m_X \Box m_Y)\circ
Tp_{TX,TY}$. Hence \eqref{eq:tau_5} follows by the third property in 
Definition \ref{def:Frob-sep} of the separable Frobenius functor $U$ and
associativity of the action $m_X\Box m_Y:T(TX\Box TY) \to TX \Box TY$.
\endproof

\begin{lemma}\label{lem:morePi}
Let $(T,\tau,\tau_0)$ be an opmonoidal endofunctor of a monoidal
category $(\M,\ox,K)$, and $(T,m,u)$ a monad on $\M$, and suppose that
equation~\eqref{eq:tau_1} holds. Then the morphism
$$\sqcap:=\big(
\xymatrix{
TK \ar[r]^-{u_{TK}}&
T^2K \ar[r]^-{\tau_{K,TK}}&
TK \ox T^2K\ar[r]^-{TK \ox m_K}&
TK \ox TK \ar[r]^-{TK \ox \tau_0}&
TK}
\big)$$
is idempotent, and the diagram 
\begin{equation}\label{eq:Pi_associ}
\xymatrix{
T^2 K \ar[r]^{T\sqcap} \ar[d]_{m_K} & T^2K \ar[r]^{m_K} & TK \ar[d]^{\sqcap} \\
TK \ar[rr]_{\sqcap} && TK }
\end{equation}
commutes.
\end{lemma}

\proof
In the diagram 
$$\xymatrix @C3pc {
TK \ar[r]^{\sqcap} \ar[d]_{u_{TK}} & 
TK \ar[rr]^{u_{TK}} && T^2K \ar[d]^{\tau_{K,TK}} \\
T^2 K \ar[r]^{\tau_{K,TK}} \ar[d]_{\tau_{K,TK}} & TK\ot T^2 K \ar[rr]^{TK\ot
  T\sqcap} \ar[d]^{TK\ot \tau_{K,TK}} & & 
TK\ot T^2 K \ar[d]^{TK\ot m_K} \\
TK\ot T^2 K \ar[r]^-{\tau_{K,K}\ot T^2 K} \ar@{=}[dr] & 
TK\ot TK\ot T^2 K \ar[r]^{TK\ot TK\ot m_K} \ar[d]^{TK\ot\tau_0\ot T^2 K} &
TK\ot TK\ot TK \ar[r]^-{TK\ot TK\ot\tau_0} \ar[d]^{TK\ot \tau_0\ot T K} & 
TK\ot TK \ar[d]^{TK\ot \tau_0} \\
& TK\ot T^2 K \ar[r]_{TK\ot m_K} & TK\ot TK \ar[r]_{TK\ot\tau_0} & TK
}$$
the squares at the bottom and the region at the top commute by naturality,
the triangle and the square above it commute since $(T,\tau,\tau_0)$ is 
opmonoidal, and the remaining region is seen to commute by taking $X=K$
in equation~\eqref{eq:tau_1} and then tensoring on the left by $TK$.  The 
composite of the top path is $\sqcap\sqcap$, and that of the bottom path 
is $\sqcap$.

As for commutativity of the displayed diagram, in the following diagram
$$\xymatrix{
TK \ar[d]_{u_{TK}} & T^2K \ar[l]_{m_K} \ar[d]^{u_{T^2 K}} \ar[r]^{T\sqcap} & 
T^2 K \ar[dr]^{m_K} \\
T^2 K \ar[d]_{\tau_{K,TK}} & T^3 K \ar[d]^{\tau_{K,T^2 K}} && TK \ar[d]^{u_{TK}} \\
TK\ot T^2 K \ar@{=}[ddd] & TK\ot T^3 K \ar[l]_{TK\ot Tm_K} \ar[d]^{TK\ot
  m_{TK}} \ar[r]^{TK\ot T^2\sqcap} & 
TK\ot T^3 K \ar[d]^{TK\ot m_{TK}} \ar[dr]^{TK\ot Tm_{K}} & T^2 K
\ar[d]^{\tau_{K,TK}} \\ 
& TK\ot T^2 K 
\ar@{=}@/_5pc/[dd]
\ar[d]^{TK\ot\tau_{K,TK}} \ar[r]^{TK\ot
  T\sqcap} & TK\ot T^2 K \ar[dr]^{TK\ot m_K} &  
TK\ot T^2 K \ar[d]^{TK\ot m_K} \\
& TK\ot TK\ot T^2 K \ar[r]^{TK\ot TK\ot m_K} \ar[d]^{TK\ot\tau_0\ot T^2 K} & 
TK\ot TK\ot TK \ar[r]^{TK\ot TK\ot\tau_0} \ar[d]^{TK\ot\tau_0\ot TK} &
TK\ot TK \ar[d]^{TK\ot\tau_0} \\
TK\ot T^2K \ar@/_3pc/[rr]_{TK\ot m_K} & TK\ot T^2 K \ar[r]^{TK\ot m_K} & TK\ot
TK \ar[r]^{TK\ot\tau_0} & TK  
}$$
the large regions at the top commute by naturality, the pentagonal region in
the middle commutes by the case $X=K$ of equation~\eqref{eq:tau_1}, and
remaining regions commute by naturality, associativity of $m$, and the
opmonoidal functor axioms.
\endproof

\begin{lemma}\label{lem:Pi}
Consider a weak bimonad $(T,m,u)$ on a monoidal category $(\M,\ox,K)$, with
opmonoidal structure $\tau_0$ in \eqref{eq:tau_0} and $\tau$ in
\eqref{eq:tau_XY}. Then the idempotent morphism  
$$
\sqcap:=\big(
\xymatrix{
TK \ar[r]^-{u_{TK}}&
T^2K \ar[r]^-{\tau_{K,TK}}&
TK \ox T^2K\ar[r]^-{TK \ox m_K}&
TK \ox TK \ar[r]^-{TK \ox \tau_0}&
TK}
\big)
$$
factorizes through an epimorphism $TK \to R$ and a section of it,
where $R$ denotes the object in $\M$ underlying the monoidal unit $(R,r)$ of
$\M^T$. 
\end{lemma}

\begin{proof}
The desired epimorphism is constructed as 
\begin{equation}
\label{eq:wbm_P}
P:=\big(
\xymatrix{
TK \ar[r]^-{Tp_0}&
TR \ar[r]^-{r}&
R}
\big)
\end{equation}
with a section 
\begin{equation}
\label{eq:wbm_I}
I:=\big(
\xymatrix{
R \ar[r]^-{u_R}&
TR \ar[r]^-{\tau_{K,R}}&
TK \ox TR \ar[r]^-{TK \ox r}&
TK \ox R \ar[r]^-{TK \ox i_0}&
TK
}
\big),
\end{equation}
where 
$(p,p_0)$ denotes the monoidal structure and 
$(i,i_0)$ denotes the opmonoidal structure
of the forgetful functor $U:\M^T \to \M$.
\end{proof}

\begin{lemma}\label{lem:E_split}
For a weak bimonad $T$ and any $T$-algebras $(A,a)$ and $(B,b)$, the
idempotent morphism
$$
E_{A,B}:=\big(
\xymatrix{
A\ox B \ar[r]^-{u_{A\ox B}}&
T(A\ox B) \ar[r]^-{\tau_{A,B}}&
TA \ox TB \ar[r]^-{a\ox b}&
A \ox B
}\big)
$$
is equal to $i_{A,B}\circ p_{A,B}$, where $\tau$ is the
natural transformation \eqref{eq:tau_XY}, $(p,p_0)$
denotes the monoidal structure and $(i,i_0)$ denotes the
opmonoidal structure of the forgetful functor $U:\M^T \to \M$.
\end{lemma}

\begin{proof}
This is immediate by \eqref{eq:tau_XY} and unitality of the $T$-actions $a$,
$b$ and $m_A \Box m_B$.
\end{proof}

The 2-category $\Mnd(\Cat)$ of monads was extended in
\cite{Bohm:weak_th_mon} to a 2-category $\Mnd^i(\Cat)$, as follows. 
The objects of $\Mnd^i(\Cat)$ are the monads in $\Cat$. 
The 1-cells from a monad $(T,m,u)$ on a category $\C$ to a monad $(T',m',u')$
on ${\mathcal C}'$,  are pairs consisting of a functor $V: \C \to \C'$ and a
natural transformation $\psi: T'V \to VT$, such that the diagram below on the
left commutes, while the 2-cells $(V,\psi) \to (W,\phi)$ are natural
transformations $\omega:V \to W$ such that the diagram on the right commutes.
\begin{equation*}
\xymatrix{
T'T'V \ar[r]^{T'\psi} \ar[d]_{m'V} & T'VT \ar[r]^{\psi T} & VTT\ar[d]^{Vm} \\
T'V \ar[rr]_{\psi} && VT 
}\quad
\xymatrix{
T'T'V \ar[r]^{T'\psi} & T'VT \ar[r]^{T'\omega T} & T'WT \ar[r]^{\phi T} & WTT
\ar[d]^{Wm} \\ 
T'V \ar[u]^{T'u'V} \ar[r]_{\psi} & VT \ar[rr]_{\omega T} && WT }
\end{equation*}

There is a variant, $\Mnd^p(\Cat)$, of this 2-category which has the same 
objects and 1-cells but in which a 2-cell $(V,\psi)\to(W,\phi)$ is a
natural transformation $\omega:V\to W$ such that  
$
\phi \circ T'\omega = Wm \circ \phi T \circ u' W T \circ \omega T \circ \psi
$.

For a monad $T$ on $\C$ and a monad $T'$ on $\C'$, we say that a functor
${\overline V}: \C^T \to \C^{\prime T'}$ is a {\em weak lifting} of a functor
$V: \C \to {\mathcal C}'$ if there exists a split natural monomorphism $i:U'
{\overline V} \to VU$ (where $U:\C^T \to \C$ and $U':\C^{\prime
  T'} \to \C'$ are the forgetful functors).
Associated to a 1-cell $(V,\psi)$ in $\Mnd^i(\Cat)$, there is an
idempotent natural transformation $VU \to VU$. Evaluated on a $T$-algebra
$(A,a)$, it is the morphism $Va\circ \psi A \circ u' VA:VA \to VA$. 
Whenever it splits, (that is, it factorizes through some natural epimorphism $VU
\to V_0$ and a section), the resulting functor $V_0:\C^T \to \C'$ has a
lifting to a functor ${\overline V}: \C^T \to \C^{\prime T'}$, which is
clearly a weak lifting of $V$.  
Conversely, every weak lifting ${\overline V}: \C^T \to \C^{\prime T'}$ of a
functor $V:\C \to\C'$ arises in this way from a unique 1-cell $(V,\psi)$ in 
$\Mnd^i(\Cat)$ such that the corresponding idempotent natural
transformation splits: see \cite[Theorem 4.4]{Bohm:weak_th_mon}.

A natural transformation ${\overline \omega}: {\overline V} \to {\overline W}$
between weakly lifted functors is said to be a {\em weak i-lifting} of a
natural transformation $\omega:V \to W$ provided that $\omega U \circ i =
i\circ U' {\overline \omega}$. 
By \cite[Proposition 4.3]{Bohm:weak_th_mon}, a natural transformation has a weak
i-lifting if and only if it is a 2-cell in $\Mnd^i(\Cat)$. Symmetrically, 
${\overline \omega}$ is said to be a {\em weak p-lifting} of $\omega$
provided that $p\circ \omega U = U' {\overline \omega}\circ p$, in terms of a
natural retraction $p$ of $i$. By \cite[Proposition 4.3]{Bohm:weak_th_mon}, a
natural transformation has a weak p-lifting if and only if it is a 2-cell in
$\Mnd^p(\Cat)$.  

\begin{theorem}\label{thm:other}
Let $(T,m,u)$ be a monad on a monoidal category $(\M,\ot,K)$, where
$(T,\tau,\tau_0):(\M,\ot,K)\to(\M,\ot,K)$ is an opmonoidal functor. 
Then equations~\eqref{eq:tau_1}--\eqref{eq:tau_5} hold if and only if
the following conditions are satisfied:
\begin{itemize}
\item[{(i)}] The functor 
$1 \stackrel{K}{\to} \M \stackrel{T}{\to} \M
$
and the natural transformation 
$T^2K \stackrel{m_K}{\to} TK \stackrel{\sqcap}{\to} TK $
constitute a 1-cell $1 \to T$ in $\Mnd^i(\Cat)$.
\item[{(ii)}] The functor 
$\M \x \M  \stackrel{\ox}{\to} \M$
and the natural transformation 
$T(\bullet \ox \bullet) \stackrel{\tau}{\to} T(\bullet) \ox
T(\bullet)$
constitute a 1-cell $T \x T \to T$ in $\Mnd^i(\Cat)$.
\item[{(iii)}] The natural transformations

$$\xymatrix @C1pc {
\M\t\M \ar[rr]^{\M\t T}  \ar@{}[rrd]|{{\displaystyle \Downarrow}\, \M\t\tau_0}
&& \M\t\M \ar[d]^{\ot} &  
\M\t\M \ar[rr]^{T\t \M}  \ar@{}[rrd]|{{\displaystyle \Downarrow}\, \tau_0\t\M}
&& \M\t\M \ar[d]^{\ot} \\  
\M \ar[u]^{\M\t K} \ar[rr]_{\M} & {} & \M & 
\M \ar[u]^{K\t\M} \ar[rr]_{\M} & {} & \M }$$

are 2-cells in $\Mnd^i(\Cat)$.
\item[{(iv)}] 
The idempotent natural transformations $E_{TX,TY}$ and $E^{(3)}_{TX,TY,TZ}$ 
which are respectively the composites
$$\xymatrix{
TX\ot TY \ar[r]^{u_{TX\ot TY}} & T(TX\ot TY) \ar[r]^{\tau_{TX,TY}} & 
T^2 X\ot T^2 Y \ar[r]^{m_X\ot m_Y} & TX\ot TY }\quad\textit{and}$$
$$\xymatrix @C3pc {
TX\ot TY\ot TZ \ar[r]^{u_{TX\ot TY\ot TZ}} & T(TX\ot TY\ot TZ)
\ar[r]^{\tau^{(3)}_{TX,TY,TZ}} & 
T^2X\ot T^2Y\ot T^2Z  
\ar[r]^{m_X\ot m_Y\ot m_Z} & TX\ot TY\ot TZ},$$
make the following diagram commute:
$$\xymatrix @R3pc @C3pc {
TX\ot TY\ot TZ \ar[r]^{E_{TX,TY}\ot TZ} \ar[dr]^{E^{(3)}_{TX,TY,TZ}}
\ar[d]_{TX\ot E_{TY,TZ}} & TX\ot TY\ot TZ \ar[d]^{TX\ot E_{TY,TZ}} \\ 
TX\ot TY\ot TZ \ar[r]_-{E_{TX,TY}\ot TZ} & TX\ot TY\ot TZ }$$
for any objects $X,Y,Z$ in $\M$.
\end{itemize}
\end{theorem}

\proof
Using the definition of 1-cells in
$\Mnd^i(\Cat)$, assertion (i) is seen to be equivalent to the equation
$\sqcap\circ m_K\circ T\sqcap = \sqcap\circ m_K$ 
which by Lemma~\ref{lem:morePi} holds whenever equation~\eqref{eq:tau_1} does.
Assertion (ii) is clearly equivalent to \eqref{eq:tau_5}. Condition (iii)
depends on the 1-cells in $\Mnd^i(\Cat)$ constructed in parts (i) and (ii). Now 
$\tau_0\circ\sqcap=\tau_0$ by opmonoidality of $T$, and then the 
two conditions in (iii) are equivalent to \eqref{eq:tau_1} and \eqref{eq:tau_2}.
Similarly the two conditions in  (iv) are equivalent to 
\eqref{eq:tau_3} and \eqref{eq:tau_4}.
\endproof

Our next aim is to prove the other half of Theorem~\ref{thm:main}, which 
we state as:

\begin{proposition}\label{prop:tau>wbm}
Consider a monad $(T,m,u)$ on a monoidal category $({\mathcal 
M},\ox,K)$ and an opmonoidal structure $(\tau,\tau_0)$ on the functor $T$. 
Assume that the identities \eqref{eq:tau_1}-\eqref{eq:tau_5}
hold and that the following idempotent morphisms split:
\begin{equation}
\label{eq:Pi}
\sqcap:=\big(
\xymatrix{
TK \ar[r]^-{u_{TK}}&
T^2K \ar[r]^-{\tau_{K,TK}}&
TK \ox T^2K\ar[r]^-{TK \ox m_K}&
TK \ox TK \ar[r]^-{TK \ox \tau_0}&
TK}
\big)
\end{equation}
and 
\begin{equation}
\label{eq:E_AB}
E_{A,B}:=\big(
\xymatrix{
A \ox B \ar[r]^-{u_{A\ox B}}&
T(A\ox B) \ar[r]^-{\tau_{A,B}}&
TA \ox TB \ar[r]^-{a\ox b}&
A \ox B
}\big),
\end{equation}
for any $T$-algebras $(A,a)$ and $(B,b)$. Then $T$ is a weak
bimonad.   
\end{proposition}

\begin{proof}
We prove this claim by constructing a monoidal structure on $\M^T$
weakly lifting that of $\M$, and by showing that with respect to this monoidal
structure the forgetful functor $U:\M^T \to \M$ is separable Frobenius.

By Theorem \ref{thm:other} (ii), $(\ox,\tau)$ is a 1-cell in $\Mnd^i(\Cat)$
(equivalently, in $\Mnd^p(\Cat)$) and so induces a weak lifting
$\bx:\M^T\x\M^T\to\M^T$; that is, a functor $\bx$ equipped with natural
transformations 
$$\xymatrix @R1pc @C1pc {
\M^T\x\M^T \ar[rr]^-{\bx} \ar[dd]_{U\x U} && \M^T \ar[dd]^{U} \\
\ar@{}[rr]|(0.6){\Downarrow i} \ar@{}[rr]|(0.4){p \Uparrow } && \\
\M\x\M \ar[rr]_-{\ox} && \M }$$   
such that $p\circ i$ is the identity natural transformation. Explicitly, for
$T$-algebras $(A,a)$ and $(B,b)$, the to-be-tensor product $(A,a)\bx
(B,b)=(A\Box B,a\Box b)$ is given by splitting the idempotent \eqref{eq:E_AB}
to obtain $A\Box B$, via maps $i_{A,B}:A\Box B\to A\ox B$ and $p_{A,B}:A\ox
B\to A\Box B$, and then $a\Box b$ is the composite
$$\xymatrix{
T(A\Box B) \ar[r]^-{Ti_{A,B}} & T(A\ox B) \ar[r]^-{\tau_{A,B}} & 
TA\ox TB \ar[r]^-{a\ox b} & A\ox B \ar[r]^-{p_{A,B}} & A\Box B. }$$

By coassociativity of $\tau$, the associativity isomorphism in $\M$
is an invertible 2-cell both in $\Mnd^i(\Cat)$ and $\Mnd^p(\Cat)$. So it
weakly lifts to an associativity isomorphism for $\bx$ such that the
following diagrams, with the associativity isomorphisms on the vertical
arrows, commute.
$$
\xymatrix@C=1.1pc{
(A \Box B)\Box C \ar[r]^-{i_{A \Box B,C}}\ar[d]&
(A\Box B) \ox C \ar[r]^-{i_{A,B}\ox C}&
(A \ox B) \ox C \ar[d]\\
A\Box (B\Box C) \ar[r]_-{i_{A,B\Box C}}&
A \ox (B \Box C) \ar[r]_-{A \ox i_{B,C}}&
A \ox (B \ox C)
}\quad \quad
\xymatrix@C=1.1pc{
(A \ox B) \ox C \ar[r]^-{p_{A,B}\ox C}\ar[d]&
(A\Box B) \ox C \ar[r]^-{p_{A \Box B,C}}&
(A \Box B)\Box C \ar[d]\\
A \ox (B \ox C) \ar[r]_-{A \ox p_{B,C}}&
A \ox (B \Box C) \ar[r]_-{p_{A,B\Box C}}&
A\Box (B\Box C)
}
$$
That is, $p$ and $i$ satisfy the associativity 
conditions that will be needed to make $U$ a monoidal and an opmonoidal
functor. The pentagon identity for $\bx$ follows from that for $\ox$ and
commutativity of either of the diagrams above.

Next we construct the unit object for the monoidal category $\M^T$. By
Theorem \ref{thm:other} (i), $(TK,\sqcap\circ m_K)$ is a 1-cell in
$\Mnd^i(\Cat)$ and so gives an object $(R,r)$ of $\M^T$. 
Explicitly, $R$ is obtained by splitting the idempotent $\sqcap$ via maps
$I:R\to TK$ and $P:TK\to R$, and $r$ is  the composite
$$\xymatrix{
TR \ar[r]^-{TI} & T^2K \ar[r]^-{m_K} & TK \ar[r]^-{P} & R. }$$

The unit constraints for $\M^T$ are constructed by applying Theorem
\ref{thm:other} (iii). The 2-cells of $\Mnd^i(\Cat)$ therein
induce morphisms $\runit_A:A\Box R\to A$ and $\lunit_A:R\Box A\to A$ of
$T$-algebras, natural in the $T$-algebra $(A,a)$. Explicitly,  $\runit_A$ and
$\lunit_A$  are given by the composites 
\begin{equation}\label{eq:lambda-rho}
\xymatrix @R1pc {
A\Box R \ar[r]^-{i_{A,R}} & A\ox R \ar[r]^-{A\ox I} & A\ox TK
\ar[r]^-{A\ox\tau_0} & A\\ 
R\Box A \ar[r]^-{i_{R,A}} & R\ox A \ar[r]^-{I\ox A} & TK\ox A \ar[r]^-{\tau_0\ox
  A} & A,}  
\end{equation}
respectively.
Since $A\Box R$ was constructed by splitting the idempotent $E_{A,R}$, to show
that $\runit_A$ is invertible, it will suffice to show that 
$$\xymatrix{
A\ox R \ar[r]^-{E_{A,R}}&A \ox R \ar[r]^-{A\ox I} & A\ox TK
\ar[r]^-{A\ox\tau_0} & A }$$  
is the epimorphism part of a splitting for $E_{A,R}$. We claim that the other
half of the splitting can be taken to be  
$$\xymatrix{
A \ar[r]^-{u_A} & TA \ar[r]^-{\tau_{A,K}} & TA\ox TK \ar[r]^-{a\ox P} & A\ox
R.}$$ 
By commutativity of the following diagram, one composite yields the identity
morphism on $A$.
$$
\xymatrix@=1.4pc{
A \ar[r]^-{u_A} \ar[rdddd]_-{u_A}
\ar@{=}@/_8em/[rrrddddd]&
TA \ar[r]^-{\tau_{A,K}} \ar[d]^-{u_{TA}} &
TA\ox TK \ar[r]^-{a\ox P}& 
A\ox R \ar[d]^-{u_{A\ox R}} \\
&T^2A\ar[r]^-{T\tau_{A,K}}\ar[ddd]^-{m_A}&
T(TA\ox TK)\ar[r]^-{T(a\ox P)}\ar[d]_-{\tau_{TA,TK}}& 
T(A\ox R) \ar[d]^-{\tau_{A,R}} \\
&&
T^2A \ox T^2K\ar[r]^-{Ta\ox TP}\ar[dd]_-{m_A \ox m_K}& 
TA \ox TR \ar[d]^-{a\ox r}\\
&&&A \ox R \ar[d]^-{A \ox I}\\
&TA\ar[r]^-{\tau_{A,K}}&
TA \ox TK\ar[r]^-{a \ox \sqcap}\ar[rd]_-{a\ox \tau_0}&
A \ox TK\ar[d]^-{A \ox \tau_0}\\ 
&&&A}
$$
The top two squares on the right hand side commute by naturality. The third
(pentagonal) region below them commutes by the associativity of $a$,
definition of $r$ and \eqref{eq:Pi_associ}. The other pentagonal region on its
left commutes by \eqref{eq:tau_5} and the region on the left of that commutes
by the unitality condition on a monad. The triangle at the bottom commutes by
the (straightforward) fact that $\tau_0\circ \sqcap=\tau_0$. The bottom-left
region commutes by the opmonoidality of $T$ and unitality of $a$. 

Composition in the opposite order yields $E_{A,R}$ by commutativity of
$$
\xymatrix@R=1pc@C=.02pc{
&A\ox R \ar@{=}[rrrrrrr]\ar[d]_-{A\ox I}&&&&&&&
A\ox R\ar[lddddd]^-{E_{A,R}}\\
&A\ox TK \ar[urrrrrrr]^-(0.2){A\ox P}\ar[rrrrrr]_-{E_{A,TK}}\ar[ddd]_-{E_{A,TK}}
\ar@(l,u)[ld]_-{u_{A \ox TK}}
&&&&&&
A\ox TK \ar[dddd]^(.2){A \ox P}
\ar[lldddd]_-{A\ox {\overline \sqcap}}
&\\
T(A \ox TK) \ar[d]_-{\tau_{A,TK}}&&&&&&&&\\
TA \ox T^2K\ar@(d,l)[rd]_-{a\ox m_K}&&&&&
&&&\\
&A \ox TK \ar[d]_-{A \ox \tau_0}&&&&&&&\\
&A\ar[rr]^-{A\ox u_K}\ar@(d,l)[rd]_-{u_A}&&
A\ox TK \ar[rr]^-{E_{A,TK}}&&
A\ox TK \ar[rr]^-{A\ox P}&&
A\ox R &\\
&&TA\ar[rr]_-{\tau_{A,K}}&&
TA \ox TK \ar@(r,d)[ur]_-{a\ox TK}&&&&
}
$$
The leftmost vertical path is equal to $(A \ox \tau_0)\circ (A\ox I)\circ
E_{A,R}$ since by the definition of $r$, $\tau_0\circ I \circ r=\tau_0 \circ
m_K \circ TI$. The regions surrounded by the curved arrows commute by
\eqref{eq:E_AB}. The triangle at the top commutes by the definitions of $I$
and $P$. The concave quadrangle below it commutes by naturality of $E$, since
$P$ is a morphism of $T$-algebras by \eqref{eq:Pi_associ}. The large polygon
at the bottom-left involves the morphism  
$$
{\overline
  \sqcap}:=
(\tau_0 \ox TK)\circ E_{TK,TK}\circ (TK\ox u_K)=
(\tau_0\ox TK)\circ (m_K\ox TK)\circ \tau_{TK,K}\circ u_{TK},
$$
where the last equality follows by \eqref{eq:E_AB}.
In order to see that this polygon commutes, note that associativity of $a$
together with \eqref{eq:tau_1} implies 
$$
E_{A,TK}\circ (A\ox \tau_0\ox TK)\circ
(E_{A,TK}\ox TK)=(A\ox \tau_0\ox TK)\circ E^{(3)}_{A,TK,TK}.
$$
Using the first one of the equivalent forms of ${\overline \sqcap}$ above, 
this implies commutativity of the bottom-left polygon. In order to see that
the triangle on its right commutes, use the second form of 
${\overline \sqcap}$. By \eqref{eq:tau_2} and opmonoidality of $T$, it obeys
$\tau_0 \circ m_K\circ T{\overline  \sqcap}= \tau_0 \circ m_K$ which implies
$\sqcap\circ {\overline \sqcap}=\sqcap$ hence commutativity of the triangle in 
question.  

The case of $\lunit$ is similar. We record here the explicit forms of 
$\runit^{-1}_A$ and $\lunit^{-1}_A$ as 
\begin{equation}\label{eq:lambda-rho-inverse}
\xymatrix @R1pc {
A \ar[r]^-{u_A} & TA \ar[r]^-{\tau_{A,K}} & TA\ox TK \ar[r]^-{a\ox P} & A\ox R
\ar[r]^-{p_{A,R}} &  
A\Box R \\
A \ar[r]^-{u_A} & TA \ar[r]^-{\tau_{K,A}} & TK\ox TA \ar[r]^-{P\ox a} & R\ox A
\ar[r]^-{p_{R,A}} & 
R\Box A, }  
\end{equation}
respectively.

To conclude that $\M^T$ is a monoidal category, we only need to prove that the
triangle condition holds. This follows by functoriality of weak lifting
because both morphisms  
$$
\xymatrix{
(A\Box R) \Box B \ar[r]&
A\Box (R \Box B) \ar[r]^-{A \Box \lunit_B}&
A \Box B
}
\quad \textrm{and}\quad
\xymatrix{
(A\Box R) \Box B \ar[r]^-{\runit_A\Box B}&
A \Box B
}
$$
are weak i-liftings of $(A \ox \tau_0)\ox B:(A \ox TK) \ox B \to A \ox B$, for
any $T$-algebras $A$ and $B$.  

It remains to show that the forgetful functor $U:\M^T\to\M$ is separable
Frobenius. We already have the binary parts of the monoidal and opmonoidal
structures, in the form of morphisms $p_{A,B}:A\ox B\to A\Box B$ and
$i_{A,B}:A\Box B\to A\ox B$. We already proved that they satisfy the
associativity, respectively, coassociativity conditions. A counit $i_0$ for
$U$, so that $(U,i,i_0)$ becomes an opmonoidal functor, is constructed as the
composite 
$$\xymatrix{
R \ar[r]^{I} & TK \ar[r]^{\tau_0} & K}$$
and the counit laws then reduce to the equations~\eqref{eq:lambda-rho} defining
$\lunit$ and $\runit$. 
The unit $p_0$ for the monoidal structure of $U$ will be the composite
$$\xymatrix{K \ar[r]^{u_K} & TK \ar[r]^{P} & R\ .}$$
One of the unit laws follows by commutativity of the diagram below; the other
is similar and left to the reader. 
$$\xymatrix{
A \ar[r]^-{A\ox u_K} \ar[d]_-{u_A} & A\ox TK \ar[r]^-{A\ox P} \ar[d]_-{u_{A\ox
    TK}} & 
A\ox R \ar[r]^-{p_{A,R}} \ar[d]_-{u_{A\ox R}} & A\Box R \ar[dd]^-{i_{A,R}} \\
TA \ar[r]^-{T(A\ox u_K)} \ar[d]_-{\tau_{A,K}} & 
T(A\ox TK) \ar[r]^-{T(A\ox P)} \ar[d]_-{\tau_{A,TK}} & 
T(A\ox R) \ar[d]_-{\tau_{A,R}} && \\
TA\ox TK \ar[r]^-{TA\ox Tu_K} \ar@{=}[dr] & 
TA\ox T^2K \ar[r]^-{TA\ox TP} \ar[d]^-{TA\ox m_K} & 
TA\ox TR \ar[r]^-{a\ox r} & A\ox R \ar[d]^-{A\ox I} \\
& TA\ox TK \ar[rr]^-{a\ox \sqcap} \ar[d]_-{TA\ox \tau_0} && 
A\ox TK \ar[d]^-{A\ox \tau_0} \\ 
& TA \ar[rr]_{a} && A}$$
The four squares in the top left corner commute by naturality; the large region 
in the top right corner by definition of $i$ and $p$; the triangle by one of
the unit laws for a monad, the region to its right by definition of $r$ and
\eqref{eq:Pi_associ}, and the bottom  region by the equation $\tau_0\circ
\sqcap=\tau_0$ once again.  The left/bottom path yields an identity morphism
by opmonoidality of $T$ and unitality of $a$.

The separability condition $p_{A,B}\circ i_{A,B}= A \Box B$ holds by
construction.   
As for the Frobenius conditions in Definition \ref{def:Frob-sep}, by 
\eqref{eq:E_AB} we have $E_{A,B\Box C}=(A\ox p_{B,C})\circ
E^{(3)}_{A,B,C}\circ(A\ox i_{B,C})$, and now the first Frobenius condition
follows by Theorem~\ref{thm:other}(iv); the other Frobenius condition is
proved similarly. 
\end{proof} 


\section{Weak bimonads vs bimonads over a separable Frobenius base}
\label{sec:Frob_sep_base}

The aim of this section is to study the category of weak bimonads on a given
Cauchy complete monoidal category $\M$. As a main result, we prove that it is 
 equivalent  to an appropriate category of bimonads on bimodule
categories over separable Frobenius monoids in $\M$. 

Recall that if T is a bimonad then $\M^T$ can be given a monoidal structure so
that the forgetful functor $U:\M^T \to \M$ is strict monoidal. Conversely, any
monad $T$ for which $U:\M^T\to \M$ is strong monoidal can be made into a
bimonad; and these two processes are, in a suitable sense, mutually
inverse. Similarly, if $g:T\to T'$ is a morphism of bimonads and $\M^T$ and
$\M^{T'}$ are made monoidal as above, then the induced functor $g^*:\M^{T'}
\to \M^T$ is strict monoidal; and conversely if $g:T \to T'$ is a morphism
of monads for which the induced functor $g^*$ is opmonoidal, compatibly with
the forgetful functors, then $g$ can be made into a morphism of bimonads.
                                         
For the entire section, we introduce the following notation. We work in a
monoidal category $\M$, with monoidal product $\ox$ and
monoidal unit $K$. For a weak bimonad $T$, the monad structure is denoted by
$m:T^2 \to T$ and $u:\M\to T$. The opmonoidal structure of $T$ is denoted by
$\tau_{X,Y}:T(X \ox Y)\to TX\ox TY$ and $\tau_0:TK \to
K$. The forgetful functor $\M^T\to \M$ is called $U$. The monoidal unit of
$\M^T$ is denoted by $(R,r)$. 
By the separable Frobenius property of $U$, $R$ is a separable Frobenius
monoid in $\M$. Its monoid structure is denoted by $(\mu:R\ox R \to
R,\eta:K\to R)$ and for the comonoid structure we write $(\delta:R \to R \ox
R,\varepsilon: R\to K)$.   
(See their explicit expressions in terms of $(m,u)$ and $(\tau,\tau_0)$
below.) 
For the monoidal category of $R$-bimodules, the forgetful functor is denoted
by $V:{}_R\M_R \to \M$. 
We use the notation $\sqcap$ introduced in \eqref{eq:Pi}, $E$ in
\eqref{eq:E_AB} and $E^{(3)}$ in 
Theorem \ref{thm:other} (iv). 
For other (weak) bimonads $T'$, ${\widetilde T}$, etc, we use the
same symbols introduced for $T$, distinguished by prime, tilde, etc. 

Our starting point is the following result due to Szlach\'anyi. 

\begin{theorem}\label{thm:strong_factor}
\cite[Theorem 2.2 and Lemma 6.2]{Szl:Bru}
Any separable Frobenius functor $U$, from a monoidal category
${\mathcal N}$ with unit $R$, to a Cauchy complete monoidal category $\M$,
factorizes through the forgetful functor ${}_{UR}{\mathcal M}_{UR} \to
{\mathcal M}$ via a strong monoidal functor ${\mathcal N} \to {}_{UR}{\mathcal
  M}_{UR}$. 
\end{theorem}

In particular, for a weak bimonad $T$ on a Cauchy complete monoidal category
$\M$, the forgetful functor $U:\M^T\to \M$ factorizes through a strong
monoidal functor ${\widetilde U}$ from $\M^T$ to the bimodule category ${}_R
\M_R$ for the monoidal unit $R$ of $\M^T$ and the forgetful functor $V:{}_R
\M_R \to \M$. Explicitly, the monoid structure of $R$ comes out as
\begin{equation}\label{eq:R_monoid}
\mu:=\big(
\xymatrix{
R \ox R \ar[r]^-{E_{R,R}}&
R \ox R \ar[r]^-{R\ox I}&
R \ox TK \ar[r]^-{R \ox \tau_0}&
R
}
\big)
\qquad 
\eta:=\big(
\xymatrix{
K \ar[r]^-{u_K}&
TK \ar[r]^-{P}&
R}
\big)
\end{equation}
and its comonoid structure is given by 
\begin{equation}\label{eq:R_comonoid}
\delta:=\big(
\xymatrix{
R \ar[r]^-{R\ox u_K}&
R \ox TK \ar[r]^-{R \ox P}&
R \ox R \ar[r]^-{E_{R,R}}&
R \ox R
}
\big)
\qquad 
\varepsilon:=\big(
\xymatrix{
R \ar[r]^-{I}&
TK \ar[r]^-{\tau_0}&
K
}
\big) .
\end{equation}
By \eqref{eq:lambda-rho}, ${\widetilde U}$ takes a
$T$-algebra $(A,a)$ to the $R\ox \bullet\ox R$-algebra $(A,\varrho_A)$ with
the structure morphism 
\begin{equation}\label{eq:R_actions}
\xymatrix{
\varrho_A=\big( R \ox A \ox R \ar[r]^-{E^{(3)}_{R,A,R}}&
R \ox A \ox R \ar[r]^-{I\ox A \ox I}&
TK \ox A \ox TK \ar[r]^-{\tau_0\ox A \ox \tau_0}&
A\big),
}
\end{equation}
where
$\xymatrix@=1.1pc{TK \ar@{->>}[r]^-P&R\ \ar@{>->}[r]^-I&TK}$
denotes a chosen splitting of the idempotent morphism $\sqcap$ of
\eqref{eq:Pi}. (Recall that $R\ox \bullet\ox R$- algebras $(M,\varrho_M)$ are
in bijection with $R$-bimodules $(M,\alpha_M:M\ox R \to M,\beta_M:R \ox M \to
M)$ via the correspondence $\varrho_M= \alpha_M \circ (\beta_M \ox R)= \beta_M
\circ (R \ox \alpha_M)$.)  

Next we compare the monadicity properties of the functors in the factorization
in Theorem \ref{thm:strong_factor}.

\begin{lemma}\label{lem:split_coeq_refl}
For a separable Frobenius monoid $R$ in a Cauchy complete monoidal category
$\M$ with forgetful functor $V:{}_R\M_R\to \M$, any $V$-contractible pair
is a split coequalizer pair. (For the terminology we refer to
\cite{BarrWells:TTT}.) 
\end{lemma}

\begin{proof}
Consider a Cauchy complete category ${\mathcal C}$ and an adjunction 
$L\dashv V:{\mathcal C}\to \M$ in which the counit $n:LV\to 1$ is split by a
natural monomorphism ${\overline n}$. Under these assumptions, any
$V$-contractible pair is a split coequalizer pair. Indeed, if for some
morphisms  $\mu,\nu:M\to N$ in ${\mathcal C}$, the first diagram in  
$$
\xymatrix@=4pc{
VM \ar@<6pt>[r]^-{V\mu} \ar@<-6pt>[r]_-{V\nu} & \ar[l]|\xi VN
}
\qquad \qquad
\xymatrix{
M \ar@<6pt>[rrr]^-{\mu} \ar@<-6pt>[rrr]_-{\nu} &&& \ar[lll]|{nM\circ L\xi \circ
  {\overline n}N} N
}
$$
is a contractible pair, then so is the second one. Hence by Cauchy
completeness of ${\mathcal C}$, the coequalizer of $\mu$ and $\nu$ exists.  

We conclude by applying this observation to the adjunction 
$R\ox \bullet \ox R\dashv V:{}_R \M_R \to \M$, whose counit is given by the
$R \ox \bullet \ox R$-action $\varrho_M:R \ox M \ox R \to M$, for any object
$(M,\varrho_M)$ of ${}_R \M_R$, hence by separable Frobenius property of $R$
it is split by $(R \ox \varrho_M \ox R)\circ (\delta \circ \eta \ox M \ox
\delta \circ \eta)$.
\end{proof}

\begin{proposition}\label{prop:W_monadic}
Let $R$ be a separable Frobenius monoid in a Cauchy complete monoidal
category $\M$. Then for a Cauchy complete category ${\mathcal C}$, a functor
$W:{\mathcal C} \to {}_R\M_R$ is monadic if and only if its composite with the
forgetful functor $V:{}_R \M_R \to \M$ is monadic.
\end{proposition}

\begin{proof}
This is proved by applying Beck's theorem \cite[Theorem 3.14]{BarrWells:TTT}.

Assume first that $W$ is monadic. Then it is immediate by monadicity of $V$
that $VW$ has a left adjoint and that it is conservative. It remains to show
that the Beck condition holds. For a $VW$-split coequalizer pair
$(\alpha,\beta)$ in ${\mathcal C}$, $(W\alpha,W\beta)$ is in
particular $V$-contractible. Hence by Lemma \ref{lem:split_coeq_refl} it is a
split coequalizer pair (evidently preserved by $V$). Then by monadicity of
$W$, there exists the coequalizer of $\alpha$ and $\beta$ and it is preserved by
$W$, so also by $VW$.  

Conversely, assume that $VW$ is monadic. 
Since $VW$ is conservative by assumption, so is $W$. As for the Beck
condition, a $W$-split coequalizer pair $(\alpha,\beta)$ is also a
$VW$-split coequalizer pair. Hence by monadicity of $VW$, its coequalizer
exists and it is preserved by $VW$. Since $V$ is faithful, this implies that
$W$ preserves the coequalizer of $\alpha$ and $\beta$. 
Thus we need only to check that $W$ has a left adjoint. This holds by a
standard adjoint-lifting argument \cite{BarrWells:TTT}, made
particularly simple here since the relevant coequalizers are split. 
In more detail, let $L$ be the left adjoint of $VW$,
with counit $n:LVW\to 1$ and unit $u:1\to VWL$.
Consider the mate of $V\varrho W$ for $\varrho:R \ox V(\bullet)\ox R
\to {}_R \M_R$ under the adjunction $L \dashv VW$; that is, the morphism 
 \begin{equation}\label{eq:q_X}
\lambda_X:=\bigl(
\xymatrix{
L(R\ox X\ox R) \ar[rr]^-{L(R\ox u_X\ox R)} & &
L(R\ox VWLX\ox R) \ar[r]^-{LV\rho_{WLX}} &
LVWLX \ar[r]^-{n_{LX}} &
LX}\bigr)
 \end{equation}
for any object $X$ of $\M$.
Whenever the coequalizer of 
$L\varrho_M, \lambda_M: L(R\ox M \ox R)\to LM$ exists for any object
$(M,\varrho_M)$ of ${}_R \M_R$, it defines a left adjoint for the lifting $W$
of $VW$; see \cite{BarrWells:TTT}.
By the separable Frobenius property of $R$, the morphism $\lambda_X$ is split
by the natural monomorphism $\lambda_{R\ox X \ox R} \circ L(\delta \circ \eta
\ox X \ox \delta \circ \eta)$. Since the diagram 
$$
\xymatrix{
L(R\ox M \ox R) \ar@<2pt>[rrrr]^-{L\varrho_M} \ar@<-2pt>[rrrr]_-{\lambda_M}
 &&&& LM \ar@{=}[d]\\
LM \ar@<2pt>[rrrr]^-{\lambda_M\circ L\big((R \ox \varrho_M \ox R)\circ
(\delta\circ \eta \ox M\ox \delta \circ \eta)\big)} \ar@<-2pt>[rrrr]_ -{LM}
\ar[u]^{\lambda_{R\ox M\ox R}\circ L(\delta\circ\eta\ox M\ox\delta\circ\eta)}
&&&& LM 
}
$$
is serially commutative, and the coequalizer of the bottom pair exists by
Cauchy completeness of $\C$, it follows that it is also a coequalizer of the
top pair, defining a left adjoint of $W$.  
\end{proof}

\begin{remark}\label{rem:T_vs_Ttilde}
Proposition \ref{prop:W_monadic} implies a relation between weak bimonads on a
Cauchy complete monoidal category $\M$ and bimonads on ${}_R\M_R$, for some
separable Frobenius monoid $R$ in $\M$. Namely, for a weak bimonad $T$, the
separable Frobenius forgetful functor $U:\M^T \to \M$ factorizes
through a strong monoidal functor ${\widetilde U}:\M^T\to {}_R \M_R$ for a
separable Frobenius monoid $R$, and the forgetful functor $V:{}_R \M_R \to
\M$; see Theorem \ref{thm:strong_factor}. By Proposition \ref{prop:W_monadic}, 
${\widetilde U}$ is also monadic hence together with its left adjoint
${\widetilde L}$, it induces a bimonad ${\widetilde T} :={\widetilde
U}{\widetilde L} $ on ${}_R \M_R$, whose Eilenberg-Moore category is
equivalent to $\M^T$; see \cite{Kelly:doctrinal}. Conversely, for a bimonad
${\widetilde T}$ on a bimodule category ${}_R\M_R$ over a separable Frobenius
monoid $R$, the composite $U$ of the forgetful functor ${\widetilde
  U}:({}_R\M_R)^{\widetilde T} \to {}_R \M_R$ and the forgetful functor
$V:{}_R \M_R\to \M$ is separable Frobenius; cf. Example
\ref{ex:Frob-sep}. It is also monadic by Proposition \ref{prop:W_monadic},
hence together with its left adjoint $L$, it induces a weak bimonad $T:=UL 
$ on $\M$ such that $\M^T$ is equivalent to $({}_R\M_R)^{\widetilde T}$. What
is more, by uniqueness of a left adjoint up to natural isomorphism, 
$T$ and $V{\widetilde T}(R \ox \bullet \ox R)$ differ by an opmonoidal
isomorphism of monads (or in fact they can be chosen equal). 
\end{remark}

\begin{remark}\label{rem:T_lift}
Consider a weak bimonad $T$ on a Cauchy complete monoidal category $\M$. By
(the proof of) Proposition \ref{prop:W_monadic}, the left adjoint ${\widetilde
L}$ of the strong monoidal functor ${\widetilde U}:\M^T\to {}_R \M_R$
(occurring in the factorization of the forgetful functor $U:\M^T\to \M$) is
constructed by choosing a splitting  
$
\xymatrix@=1.3pc{
LM \ar@{->>}[r]&
{\widetilde L}(M,\varrho_M) 
\ \ar@{>->}[r]&
LM
}
$
of the idempotent natural transformation
\begin{equation}\label{eq:Q_idemp}
\lambda_M \circ L\big((R \ox \varrho_M\ox R)\circ (\delta\circ \eta\ox M \ox
\delta\circ \eta)\big):LM\to LM,
\end{equation}
for any object $(M,\varrho_M)$ of ${}_R \M_R$, where $\lambda_M$ is as in 
\eqref{eq:q_X}. Applying $U$, this yields a split idempotent natural
transformation
$ULV = \xymatrix@=1.3pc{
TV \ar@{->>}[r]^-{q}&
U{\widetilde L}=V{\widetilde T} 
\ \ar@{>->}[r]^-{j}&
TV}=ULV$.
What is more, also as a monad, ${\widetilde T}={\widetilde U}{\widetilde L}$
is a weak $q$-lifting of $T$ hence by \cite[Proposition
3.7]{Bohm:weak_th_mon}, the Eilenberg-Moore categories $\M^T$ and
$({}_R\M_R)^{\widetilde T}$ are in fact isomorphic.  
Explicitly, there is an isomorphism  $\Xi:\M^T\to ({}_R\M_R)^{\widetilde
T}$, taking a $T$-algebra $(A,a)$ to the $R\ox \bullet \ox R$-algebra $A$
described  in \eqref{eq:R_actions}, 
with a ${\widetilde T}$-algebra
structure provided by the unique morphism ${\widetilde a}:{\widetilde T}A\to
A$ for which ${\widetilde a}\circ q_A =a$. On the morphisms, $\Xi$ acts as the
identity map. The inverse of $\Xi$ takes an object $((A,\varrho_A),{\widetilde
  a})$ of $({}_R\M_R)^{\widetilde T}$ to the $T$-algebra $A$, with structure
morphism $TA \1c {q_A} V{\widetilde T}A \1c{\widetilde a} A$, and it also acts
on the morphisms as an identity map. In particular, also the Eilenberg-Moore
categories $({}_R \M_R)^{\widetilde T}$ and $\M^{V{\widetilde T}(R \ox \bullet
  \ox R)}$ are isomorphic, for any bimonad ${\widetilde T}$ on a bimodule
category ${}_R \M_R$ over a separable Frobenius monoid $R$.
\end{remark}

The final aim of this section is to extend the correspondence in Remark
\ref{rem:T_vs_Ttilde} between weak bimonads on one hand, and bimonads over
separable Frobenius base monoids on the other hand, to an 
 equivalence 
of  categories.

\begin{definition}\label{def:cat_wbm}
A {\em morphism of weak bimonads} on a monoidal category $\M$ is
defined as an opmonoidal morphism of monads; that is, as a natural
transformation $g:T\to T'$ which is opmonoidal in the sense that, for any
objects $X$ and $Y$ in $\M$, 
$$
\xymatrix{
T(X \ox Y) \ar[r]^-{g_{X\ox Y}}\ar[d]_-{\tau_{X,Y}}
& T'(X \ox Y)\ar[d]^-{\tau'_{X,Y}}\\
TX \ox TY \ar[r]_-{g_X \ox g_Y} &
T'X \ox T'Y}
\qquad
\xymatrix{
TK\ar[r]^-{g_K}\ar[d]_-{\tau_0}&
T'K\ar[d]^-{\tau'_0}\\
K \ar@{=}[r] &
K
}
$$
and which is a morphism of monads in the sense that, for any object $X$ in
$\M$,   
$$
\xymatrix{
T^2 X \ar[r]^-{Tg_X}\ar[d]_-{m_X}&
TT'X \ar[r]^-{g_{T'X}}& 
T^{\prime 2}X \ar[d]^-{m'_X} \\
TX \ar[rr]_-{g_X}&&
T'X}
\qquad 
\xymatrix{
X \ar@{=}[r]\ar[d]_-{u_X}&
X \ar[d]^-{u'_X}\\
TX \ar[r]_-{g_X} &
T'X\ .}
$$
Weak bimonads on $\M$ (as objects) and their morphisms (as arrows) constitute 
a category $\wbm(\M)$, which contains the category of bimonads on $\M$ as a
full subcategory.
\end{definition}

\begin{example}\label{ex:wbm_morphism}
Any monoid $R$ in a monoidal category $\M$, induces a monad $R \ox \bullet \ox
R$ on $\M$.  
If $R$ is a separable Frobenius monoid in a Cauchy complete monoidal
category $\M$, then $R\ox \bullet\ox R$ is a weak bimonad; see Example
\ref{ex:Frob-sep} (3). Its opmonoidal structure is provided by  the 
maps 
$$\xymatrix@C=2pc{
R\ot R 
\ar[r]^-{\varepsilon\circ \mu} & K }
\quad \textrm{and}\quad
\xymatrix @C=7pc {
R\ot X\ot Y\ot R 
\ar[r]^-{R\ot X\ot \delta\circ \eta \ot Y\ot R} &
R\ot X \ot R\ot R\ot Y\ot R, }$$
for any objects $X,Y$ in $\M$, where $(\mu,\eta)$ and $(\delta,\varepsilon)$
denote the monoid and comonoid structures of $R$, respectively.
A morphism $\gamma:R \to R'$ of separable Frobenius monoids in $\M$ induces
a morphism of weak bimonads $\gamma\ox \bullet \ox \gamma: R\ox \bullet \ox
R\to R'\ox \bullet \ox R'$. 
\end{example}

\begin{lemma}\label{lem:R_morphism}
Consider any morphism $g:T \to T'$ of weak bimonads on a Cauchy complete
monoidal category $\M$, with monoidal units $R$, respectively $R'$, of the
Eilenberg-Moore categories $\M^T$ and $\M^{T'}$. 
There is a unique isomorphism $\gamma: R \to R'$ of separable Frobenius
monoids such that the following diagram of functors commutes,
\begin{equation}\label{eq:mon_functors}
\xymatrix{
({}_{R'}\M_{R'})^{\widetilde T'}\cong \M^{T'}\ar[rr]^-{g^*}
  \ar[d]_{\widetilde U'} &&  
\M^T \cong ({}_R \M_R)^{\widetilde T} 
\ar[d]^{\widetilde U} \\
{}_{R'}\M_{R'} \ar[rr]^-{\gamma^*}\ar[rd]_-{V'}&&
{}_R \M_R \ar[ld]^-V\\
& \M &
}
\end{equation}
where the bimonads ${\widetilde T}$ and ${\widetilde T}'$ are associated to
the weak bimonads $T$ and $T'$ as in Remark \ref{rem:T_vs_Ttilde}.
\end{lemma}

\begin{proof}
By Lemma \ref{lem:Pi}, for weak bimonads $T$ and $T'$, the associated
idempotent morphisms $\sqcap$ and $\sqcap'$ in \eqref{eq:Pi} split through $R$
and $R'$, respectively; thus there are epi-mono pairs 
$
\xymatrix@=1.3pc{
TK\ar@{->>}[r]^-{P}&
R \ \ar@{>->}[r]^-{I}&
TK
}
$ 
and 
$
\xymatrix@=1.3pc{
T'K \ar@{->>}[r]^-{P'}&
R'\ \ar@{>->}[r]^-{I'} &
T'K
}
$. Using the fact that a morphism $g:T\to T'$ of weak bimonads is an
opmonoidal 
natural transformation as well as a morphism of monads, one checks that for
any objects $X,Y$ in $\M$,  the two diagrams on the left 
\begin{equation}
\label{eq:g_id}
\xymatrix{
TK \ar[r]^{g_K} \ar[d]_{\sqcap} & T'K \ar[d]^{\sqcap'} \\
TK \ar[r]_{g_K} & T'K }
\qquad 
\xymatrix{
TX\ox TY \ar[r]^{g_X\ox g_Y} \ar[d]_{E_{TX,TY}} & T'X\ox T'Y
\ar[d]^{E'_{T'X,T'Y}} \\ 
TX\ox TY \ar[r]_{g_X\ox g_Y} & T'X\ox T'Y}
\qquad
\xymatrix{
TK \ar[r]^{g_K} & T'K \ar[d]^{P'} \\
R \ar[u]^{I} \ar[r]_{\gamma} & R' }
\end{equation}
commute and so the morphism $\gamma$ defined by the diagram on the
right is compatible both
with the monoid and with the comonoid structures of $R$ and $R'$, 
written out explicitly in \eqref{eq:R_monoid} and \eqref{eq:R_comonoid}. 
That is, $\gamma$ is a  morphism of separable Frobenius monoids. By
\cite[Proposition A.3]{PaSt:WH_mon}, any morphism of Frobenius monoids is an
isomorphism hence so is $\gamma$. 
It obviously renders commutative the lower triangle in
\eqref{eq:mon_functors}. It renders commutative also the upper square 
by commutativity of the following diagram, for any $T'$-algebra
$(A,a)$. 
\begin{equation}\label{eq:gamma_diagr}
\xymatrix{
R \ox A\ox R \ar[r]^-{u_{R\ox A \ox R}}\ar[d]_-{\gamma\ox A\ox \gamma}&
T(R\ox A\ox R)\ar[r]^-{\tau^{(3)}_{R,A,R}}\ar[d]^-{T(\gamma\ox A\ox
  \gamma)}&
TR\ox TA\ox TR \ar[r]^-{TR\ox g_{A} \ox TR}\ar[d]^-{T\gamma \ox TA\ox
  T\gamma}&
TR \ox T'A\ox TR\ar[d]^-{r\ox a \ox r}\\
R'\ox A\ox R'\ar[r]^-{u_{R'\ox A\ox R'}}\ar[d]_-{u'_{R'\ox A\ox R'}}&
T(R'\ox A\ox R')\ar[r]^-{\tau^{(3)}_{R',A,R'}}\ar[dl]_-{g_{R'\ox A\ox R'}}&
TR'\ox TA\ox TR'\ar[dl]_-{g_{R'}\ox g_{A} \ox g_{R'}}&
R\ox A\ox R\ar[dl]_-{\gamma\ox A\ox \gamma}\ar[d]^-{\varepsilon \ox A\ox
  \varepsilon} \\
T'(R'\ox A\ox R') \ar[r]_-{\tau^{\prime(3)}_{R',A,R'}}&
T'R'\ox T'A\ox T'R' \ar[r]_-{r'\ox a \ox r'}&
R'\ox A\ox R'\ar[r]_-{\varepsilon'\ox A\ox \varepsilon'}&
A
}
\end{equation}
The two squares in the upper left corner commute by naturality and the square
below them commutes by the opmonoidality of $g$. The triangles in the bottom
row commute since $g$ is a monad morphism, and since $\gamma$ is a comonoid
morphism, respectively. The remaining region commutes by commutativity of the
following diagram, 
\begin{equation}\label{eq:g_gamma}
\xymatrix{
&&&& \\ 
TR \ar[r]^-{TI}\ar[d]_-{TI}
\ar@{}[drrrr]|{\eqref{eq:Pi_associ}}    
\ar@(u,u)[rrrr]^-r
\ar@(l,l)[ddd]_-{T\gamma}&
TTK \ar[rr]^-{m_K}&&
TK \ar[r]^-{P} &
R \ar[d]^-{I} 
\ar@(r,r)[ddd]^-\gamma\\
TTK \ar[d]_-{Tg_K}\ar@{}[rd]|{\eqref{eq:g_id}}&
TTK \ar[lu]_-{TP} \ar[rrr]^-{m_K}\ar[d]^-{Tg_K}&&&
TK \ar[d]^-{g_K}\\
TT'K \ar[d]_-{TP'}&
TT'K \ar[ld]_-{TP'}\ar[r]^-{g_{T'K}}&
T'T'K \ar[ld]_-{T'P'}\ar[rr]^-{m'_K}\ar@{}[rrd]|{\eqref{eq:Pi_associ}}&&
T'K\ar[d]^-{P'}\\
TR' \ar[r]_-{g_{R'}}&
T'R' \ar[r]_-{T'I'} 
\ar@(d,d)[rrr]_-{r'}&
T'T'K \ar[r]_-{m'_K}&
T'K \ar[r]_-{P'}&
R' &} 
\end{equation}
where the undecorated region in the middle commutes since $g$ is a
morphism of monads. 
It remains to show that $\gamma$ is unique with the stated property. The upper
part of the diagram in \eqref{eq:mon_functors} commutes if and only if, for
any $T'$ algebra $(A,a)$ and the corresponding $T$-algebra $(A,a\circ g_A)$,
the $R\ox \bullet \ox R$-action $\varrho_A$ and the $R'\ox \bullet \ox
R'$-action $\varrho'_A$ in \eqref{eq:R_actions} obey 
$\varrho_A
=\varrho'_A\circ (\gamma \ox A \ox \gamma)$; 
see \eqref{eq:gamma_diagr}. 
Applying it in the case where $(A,a)$ is the monoidal unit
$(R',r')$ of $\M^{T'}$ and composing the resulting identity on the right with
$\eta \ox \eta' \ox R$, we deduce that
$\gamma=(R'\ox \varepsilon)\circ E_{R',R}\circ (\eta'\ox R)$. Thus
$\gamma$ renders commutative the last diagram in \eqref{eq:g_id} by the forms
of $P'$ and $I$ in \eqref{eq:wbm_P} and \eqref{eq:wbm_I}, the form of
$E_{R',R}$ in \eqref{eq:E_AB}, and naturality.  
\end{proof}

On the other hand, any isomorphism of separable Frobenius monoids clearly 
induces a strict monoidal isomorphism between the categories of bimodules,
which can be seen as a particular case of the following:

\begin{lemma}\label{lem:str_mon_functor}
If $g:T\to T'$ is a morphism of weak bimonads on a Cauchy complete monoidal
category $\M$ then the induced functor $g^*:\M^{T'}\to \M^T$ is strong
monoidal and \eqref{eq:mon_functors} is a commutative diagram of separable 
Frobenius monoidal
functors. 
\end{lemma}

\begin{proof}
If $T$ and $T'$ are bimonads then the monoidal structures on $\M^T$ and 
$\M^{T'}$ are lifted from that on \M, and so clearly $g^*:\M^{T'}\to\M^T$ is
strong (in fact strict) monoidal. If $T$ is only a weak bimonad then the
monoidal 
structure on $\M^T$ is only weakly lifted from that on \M, 
but it is lifted, up to equivalence, from that on ${}_R\M_R$. Thus if
$g:T\to T'$ is a morphism of weak bimonads, then the isomorphism $\gamma:R\to
R'$ of the previous lemma induces a strict monoidal isomorphism 
${}_{R'}\M_{R'}\to{}_R\M_R$, and now the strong monoidal structure on the 
composite $\M^{T'}\to{}_{R'}\M_{R'}\to{}_R\M_R$ lifts to a strong monoidal
structure on $g^*:\M^{T'}\to\M^T$. Explicitly, this is given by 
$$
\gamma_{A,B}:=\big(
\xymatrix{
A \Box B \ar[r]^-{i_{A,B}} &
A \ox B \ar[r]^-{p'_{A,B}} &
A \Box' B
}\big),
$$
and $\gamma:R\to R'$.
\end{proof}

We now turn to our category of bimonads on categories of bimodules over
separable Frobenius monoids in $\M$.

\begin{definition}
\label{def:fsbm}
For a monoidal category $\M$, the category $\fsbm(\M)$ is defined to have {\em
  objects} which are pairs $(R,{\widetilde T})$, consisting 
of a separable Frobenius monoid $R$ in $\M$ and a bimonad ${\widetilde T}$
on ${}_R\M_R$. {\em Morphisms} $(R,{\widetilde T})\to (R',{\widetilde T'})$
are pairs $(\gamma,\Gamma)$, consisting of an isomorphism $\gamma:R \to R'$ of
separable Frobenius monoids (inducing a strong monoidal isomorphism
$\gamma^*:{}_{R'}\M_{R'} \to {}_{R}\M_{R}$), and a morphism of bimonads
$\Gamma:{\widetilde T}\to \gamma^* {\widetilde T'} (\gamma^*)^{-1}$ (that is,
an opmonoidal morphism of monads, in the sense explicated
in Definition \ref{def:cat_wbm}). 
\end{definition}

There is an evident functor $\Psi:\fsbm(\M)\to\wbm(\M)$ defined as
follows. The object map is given by associating to a pair $(R,{\widetilde T})$ 
the weak bimonad induced by the composite of the forgetful functors
${\widetilde U}:({}_R \M_R)^{\widetilde T} \to {}_R \M_R$ and $V:{}_R \M_R \to
\M$, as described in Remark \ref{rem:T_vs_Ttilde}; 
explicitly, $\Psi(R,{\widetilde T})X=V\tilde{T}(R\ox X\ox R)$. A
morphism $(\gamma,\Gamma):(R,\tilde{T})\to(R',\tilde{T'})$ in $\fsbm(\M)$
gives rise to a commutative diagram of functors
$$
\xymatrix{
\M^{\Psi(R',{\widetilde T}')}\cong ({}_{R'}\M_{R'})^{\widetilde{T}'} 
\ar[d]_{\widetilde U'} \ar[rr]^{\Gamma^*} && 
({}_{R}\M_R)^{\widetilde{T}} \cong \M^{\Psi(R,{\widetilde T})}  
\ar[d]^{\widetilde U} \\
{}_{R'}\M_{R'} \ar[rr]^{\gamma^*} \ar[dr]_{V'} && {}_{R}\M_{R} \ar[dl]^{V}
\\
& \M }
$$
and so in particular to a morphism of
monads $g:\Psi(R,{\widetilde T})\to \Psi(R',{\widetilde T}')$,
explicitly, $g_X:V{\widetilde T}(R \ox X \ox R) \to V'{\widetilde T}'(R'\ox X
\ox R')$ is given by 
$$\xymatrix @C1pc {
V\widetilde{T}(R\ox X\ox R) \ar[rr]^-{V\Gamma_{R\ox X \ox R}} &&
V\gamma^*\widetilde{T}'(\gamma^*)^{-1}(R\ox X\ox R) \ar@{=}[d] \\
&& V'\widetilde{T}'(\gamma^*)^{-1}(R\ox X\ox R)
\ar[rrr]^-{V'\widetilde{T}'(\gamma\ox X\ox \gamma)} &&&
V'\widetilde{T}'(R'\ox X\ox R') }$$
and it is opmonoidal since $V$, $V'$, and $\widetilde{T}'$ are opmonoidal
functors, and $\Gamma$ and $\gamma\ox X\ox\gamma:(\gamma^*)^{-1}(R\ox
X\ox R)\to R'\ox X\ox R'$ are opmonoidal natural transformations.

\begin{theorem}
If \M is a Cauchy complete monoidal category, the functor 
$\Psi:\fsbm(\M)\to\wbm(\M)$ is an equivalence of categories.
\end{theorem}

\proof
First we show that $\Psi$ is fully faithful. Suppose then that objects
$(R,\widetilde{T})$ and $(R',\widetilde{T}')$ of $\fsbm(\M)$ are given.
We must show that any morphism $g:\Psi(R,{\widetilde T})\to
\Psi(R',{\widetilde T'})$ of weak bimonads is induced by a 
unique morphism $(\gamma,\Gamma):(R,\widetilde{T})\to(R',\widetilde{T}')$ in 
$\fsbm(\M)$. 
The existence of a unique isomorphism $\gamma:R \to R'$ of monoids, inducing an
isomorphism $\gamma^*:{}_{R'}\M_{R'}\to {}_R \M_R$ of categories rendering
commutative \eqref{eq:mon_functors}, is given by Lemma~\ref{lem:R_morphism}. 
By Lemma \ref{lem:str_mon_functor}, $g$ induces a
strong monoidal functor $g^*:\M^{T'}\to\M^{T}$. By commutativity of the upper
square in \eqref{eq:mon_functors} as a diagram of strong monoidal functors, it
is necessarily of the form $\Gamma^*$ for a unique opmonoidal monad morphism
$\Gamma:{\widetilde T}\to \gamma^*{\widetilde T}'(\gamma^*)^{-1}$.  
This proves that $\Psi$ is fully faithful. It is essentially
surjective on objects by Remark~\ref{rem:T_vs_Ttilde}.  
\endproof

Bimonads are monads in the 2-category $\textsf{OpMon}$ of monoidal
categories, opmonoidal functors and opmonoidal natural transformations;
cf. \cite{McCr:opmon_mon}. That is, they can be regarded as 0-cells in the
2-category $\Mnd( \textsf{OpMon})$. Clearly, for a Cauchy complete monoidal
category $\M$, the category $\fsbm(\M)$ is a subcategory in the opposite of
the category underlying $\Mnd(\textsf{OpMon})$.  
We may consider also the full subcategory of the underlying category of
$\Mnd(\textsf{OpMon})$, with objects the bimonads on bimodule categories over 
separable Frobenius monoids $R$ in Cauchy complete monoidal categories. In
this way (using the correspondence in Remark \ref{rem:T_vs_Ttilde}), we can
define more general morphisms between weak bimonads than the arrows in
$\wbm(\M)$ for a fixed $\M$. These more general morphisms do not need to
preserve the underlying separable Frobenius monoid $R$.


\section{An example: Weak bimonoids in braided monoidal categories}
\label{sec:weak_bimonoid}

In this section we show that weak bimonoids in a Cauchy complete braided
monoidal category $\M$ induce weak bimonads on $\M$.

\begin{theorem}
\label{thm:w_bimonoid}
For a monoid $(B,\mu,\eta)$ in a Cauchy complete braided monoidal category
$(\M,\ox,K,c)$, there is a bijection between 
\begin{itemize}
\item[{(1)}] weak bimonoids of the form $(B,\mu,\eta, \delta,\varepsilon)$ in
  $\M$;
\item[{(2)}] weak bimonads $(\bullet \ox B, \bullet \ox \mu, \bullet
  \ox \eta,\tau,\tau_0)$ on $\M$ for which the diagram
\begin{equation}
\label{eq:w_bimonoid_tau}
\xymatrix @C3pc {
X\ox Y \ox B \ar[r]^-{X\ox Y\ox \tau_{K,K}} \ar[dr]_{\tau_{X,Y}} & 
X\ox Y\ox B\ox B \ar[d]^{X\ox c_{Y,B}\ox B} \\
& X\ox B\ox Y\ox B
}
\end{equation}
commutes for all objects $X,Y$ of $\M$.
\end{itemize}
\end{theorem}

\begin{remark}
Consider a monoid $B$ in a Cauchy complete braided monoidal category
$(\M,\ox,$ $K,c)$ such that $\bullet \ox B$ is a weak bimonad on $\M$. By
naturality, for all morphisms $f:K \to X$, $g:K\to Y$, $h:K\to B$, the natural
transformation $\tau_{X,Y}: X \ox Y \ox B \to X \ox B \ox Y \ox B$ makes
$$
\xymatrix@C=3pc @R=1pc{
K \ar[dr]^(.7){f\ox g\ox h} \ar[d]_{h} \ar[r]^-{f\ox g\ox h} & 
X\ox Y\ox B \ar[r]^-{X\ox Y\ox \tau_{K,K}} & 
X\ox Y\ox B\ox B \ar[dd]^-{X\ox c_{Y,B}\ox B} \\
B \ar[d]_-{\tau_{K,K}} &
X\ox Y\ox B \ar[dr]^-{\tau_{X,Y}} &  \\
B\ox B \ar[rr]_-{f\ox B\ox g\ox B}&&
X\ox B\ox Y\ox B 
}
$$ 
commute.
Hence \eqref{eq:w_bimonoid_tau} holds provided that the monoidal unit $K$ is
a {\em `cubic generator'} in the following sense: If, for some morphisms
$p,q:X \ox Y \ox Z \to W$ in $\M$, the equality $p\circ (f\ox g \ox h)=q \circ
(f\ox g \ox h)$ holds, for {\em all} morphisms $f:K \to X$, $g:K\to Y$,
$h:K\to Z$, then $p=q$. 

The monoidal unit is a `cubic generator', for example, in the symmetric
monoidal category $\mathsf{Mod}(k)$ of modules over a commutative ring
$k$. With this observation in mind, Theorem \ref{thm:w_bimonoid} includes
Szlach\'anyi's description in \cite[Corollary 6.5]{Szl:Bru} of {\em weak
  bialgebras} over $k$ as {\em weak bimonads} on $\mathsf{Mod}(k)$. 
\end{remark}

\begin{proof}[Proof of Theorem \ref{thm:w_bimonoid}]
Suppose that $(B,\mu,\eta,\delta,\epsilon)$ is a weak bimonoid in \M.
By \cite[Proposition 3.8]{PaSt:WH_mon}, the category of $B$-modules is
monoidal and there is a strong monoidal functor from it to a certain bimodule
category ${}_R \M_R$. Furthermore, the resulting monoid $R$ is a
separable Frobenius monoid by \cite[Proposition 1.4]{PaSt:WH_mon}. In view
of Example \ref{ex:Frob-sep}, this proves that $\bullet \ox B$ is a weak
bimonad. Its opmonoidal structure comes out with $\tau_{X,Y}$ equal to 
the composite
$$\xymatrix@C=4pc{
X\ox Y\ox B \ar[r]^{X\ox Y\ox\delta} & X\ox Y\ox B\ox B \ar[r]^{X\ox c_{Y,B}\ox B} &
X\ox B\ox Y\ox B}$$
and $\tau_0=\epsilon$. Hence
\eqref{eq:w_bimonoid_tau} is satisfied.

Assume conversely that (2) holds. We claim that
$(B,\mu,\eta,\delta:=\tau_{K,K},\varepsilon:=\tau_0)$ is a weak bimonoid in
$\M$. The functor $\bullet\ox B$ is opmonoidal and so sends comonoids to 
comonoids; in particular, it sends the comonoid $K$ in $\M$ to a comonoid,
which turns out to be $(B,\delta,\epsilon)$.
Use \eqref{eq:w_bimonoid_tau} to write $\tau_{X,Y}$ as 
$(X \ox c_{Y,B}\ox B)\circ (X \ox Y \ox \delta)$, for any objects $X,Y$ of
$\M$. Substituting this expression in conditions
\eqref{eq:tau_1}-\eqref{eq:tau_5}, we obtain the following commutative diagrams.
\begin{equation}\label{eq:w_bimonoid_1}
\xymatrix @R=1pc @C=1pc{
B^2 \ar[rr]^{B\ox\eta\ox B} \ar[d]^{B\ox\delta}&&
B^3 \ar[rr]^{B\ox\delta\ox B} && B^4 \ar[rr]^{B\ox c^{-1}_{B,B}\ox B} && 
B^4 \ar[d]_{\mu\ox\mu}\\
B^3 \ar[d]^{B\ox c^{-1}_{B,B}}&&&&&&
B^2 \ar[d]_{\epsilon\ox B} \\ 
B^3 \ar[rrr]_{\mu\ox B} &&&
B^2 \ar[rrr]_{\epsilon\ox B} &&&
B}
\qquad
\xymatrix @R=1pc @C=1.5pc{
B^2 \ar[rr]^{B\ox\eta\ox B}\ar[dd]^{B\ox\delta} &&
B^3 \ar[rr]^{B\ox\delta\ox B} &&
B^4 \ar[d]_{\mu\ox\mu}\\
&&&&B^2 \ar[d]_-{\epsilon\ox B}\\
B^3 \ar[rr]_{\mu\ox B} &&
B^2 \ar[rr]_{\epsilon\ox B}&&
B
}
\end{equation}
\begin{equation}\label{eq:w_bimonoid_3}
\xymatrix @R=2pc @C=1pc{
\ K \ar[d]^{\eta\ox\eta} \ar[rrr]^{\eta} &&&
\ B \ar[rrr]^-{\delta}&&&
\ B^2 \ar[d]_{\delta\ox B} \\
B^2 \ar[rr]_{\delta\ox\delta}&&
B^4 \ar[rr]_{B\ox c^{-1}_{B,B}\ox B} && 
B^4 \ar[rr]_{B\ox \mu\ox B}&&
B^3}
\qquad
\xymatrix @R=2pc @C=1.5pc{
\ K \ar[d]^{\eta\ox\eta} \ar[rr]^{\eta} &&
\ B \ar[rr]^{\delta} &&
\ B^2 \ar[d]_{\delta\ox B}\\
B^2 \ar[rr]_{\delta\ox\delta}&&
B^4 \ar[rr]_{B\ox\mu\ox B}&&
B^3
}
\end{equation} 
\begin{equation}\label{eq:w_bimonoid_5}
\xymatrix@C=3pc{
B^4 \ar[rr]^{B\ox c_{B,B}\ox B} && B^4 \ar[d]_{\mu\ox\mu} \\
B^2 \ar[u]_{\delta\ox\delta} \ar[r]_{\mu} & B \ar[r]_{\delta} & B^2 } 
\end{equation}
Condition \eqref{eq:w_bimonoid_5} is identical to axiom (b), 
and the identities in \eqref{eq:w_bimonoid_3}
are identical to axiom (w) in the definition of a weak bimonoid in 
\cite[Definition 2.1]{PaSt:WH_mon}. Thus we only need to show that,
whenever the diagrams in \eqref{eq:w_bimonoid_1}, \eqref{eq:w_bimonoid_3},
and \eqref{eq:w_bimonoid_5} commute, 
then axiom (v) in \cite[Definition 2.1]{PaSt:WH_mon} holds; that is, the
following diagram commutes.
\begin{equation}\label{eq:ax_v}
\xymatrix@C=4pc @R=2pc{
B^3 \ar[r]^-{B\ox \delta \ox B}\ar[dd]_-{B \ox \delta \ox B}
\ar[rd]^-{\mu^2}&
B^4 \ar[r]^-{B \ox c_{B,B}^{-1}\ox B}&
B^4 \ar[d]^-{\mu \ox \mu}\\
&B \ar[rd]^-{\varepsilon}&
B^2 \ar[d]^-{\varepsilon\ox \varepsilon}\\
B^4 \ar[r]_-{\mu \ox \mu}&
B^2 \ar[r]_-{\varepsilon\ox\varepsilon}&
K}
\end{equation} 
Commutativity of the lower triangle in \eqref{eq:ax_v} follows by 
commutativity of 
$$
\xymatrix@C=3.5pc @R=1pc{
B^3 \ar@{=}[r] \ar[dd]_-{B\ox \delta \ox B} &
B^3 \ar[r]^-{B\ox \mu}\ar[d]^-{B \ox \eta \ox B^2}&
B^2 \ar@{=}[r]\ar[d]^-{B\ox \eta \ox B}&
B^2 \ar[r]^-{\mu}\ar[d]^-{B\ox \delta}&
B \ar[dddd]^-{\varepsilon}\\
&
B^4 \ar[d]^-{B\ox \delta \ox B^2}&
B^3\ar[d]^-{B\ox \delta \ox B}\ar@{}[rdd]|{\eqref{eq:w_bimonoid_1}}&
B^3 \ar[dd]^-{\mu \ox B}&\\
B^4\ar[dd]_-{\mu \ox B\ox B}\ar@{}[r]|{\eqref{eq:w_bimonoid_1}}&
B^5 \ar[d]^-{\mu \ox \mu \ox B}&
B^4 \ar[d]^-{\mu \ox \mu}&\\
&
B^3 \ar[d]^-{\varepsilon\ox B \ox B}&
B^2 \ar[d]^-{\varepsilon\ox B}&
B^2 \ar[d]^-{\varepsilon \ox B}&\\
B^3 \ar[r]_-{\varepsilon\ox B \ox B}&
B^2 \ar[r]_-\mu&
B \ar@{=}[r]&
B \ar[r]_-{\varepsilon}&
K.
}
$$
The undecorated regions commute by naturality, counitality of $\delta$ and by
associativity of $\mu$. 
Commutativity of the upper triangle in \eqref{eq:ax_v} is proved similarly,
making use of the first identity in \eqref{eq:w_bimonoid_1}. 
\end{proof}


\section{The antipode}
\label{sec:antipode}

The first attempt to equip Moerdijk's bimonad with an antipode, i.e. to
define a Hopf monad, was made by Brugui\`eres and Virelizier in
\cite{BV:Hopf_mon}. Here the authors studied bimonads on autonomous monoidal
categories such that the (left/right) duals lift to the Eilenberg-Moore
category. This generalizes finite dimensional Hopf algebras to the categorical
setting.  

A more general notion of Hopf monad was introduced in \cite{BLV} (see also
\cite{ChiklackStr:Homf_mon}). 
This is based on the observation of Lawvere \cite{Law} that a right adjoint
preserves internal homs precisely when Frobenius reciprocity holds;
this Frobenius reciprocity condition also appeared in
\cite[Theorem \& Definition 3.5]{Sch:x_R} in the context of Takeuchi
bialgebroids.
Based on this result, the following definition was proposed also
for monoidal categories which are not necessarily closed. 
For any monad $(T,m,u)$ on a monoidal category $(\M,\ox,K)$ such that $T$
admits an opmonoidal structure $(\tau,\tau_0)$ (hence in particular for any
(weak) bimonad $T$ on $\M$), there is a canonical natural   transformation,
given for any objects $X,Y$ of $\M$ by  
\begin{equation}
\label{eq:can}
\can_{X,Y}:=\big(
\xymatrix{
T(TX \ox Y)\ar[r]^-{\tau_{TX,Y}}&
T^2X \ox TY \ar[r]^-{m_X \ox TY} &
TX \ox TY
}\big).
\end{equation}
By the terminology in  \cite{BLV}, a bimonad is called a {\em right 
Hopf monad} whenever the associated natural transformation \eqref{eq:can} is
invertible. Similarly, a bimonad is a {\em left Hopf monad} when the analogous
natural transformation $T(X\ox TY)\to TX\ox TY$ is invertible, and a {\em Hopf
  monad} when it is both left and right Hopf. 

In this section we propose a definition of a {\em right weak Hopf monad} $T$
on a monoidal category $\M$ -- characterized by the property that, whenever
$\M$ is also Cauchy complete, the associated bimonad ${\widetilde T}$ (on
another monoidal category) in Remark \ref{rem:T_vs_Ttilde}, is a right Hopf
monad, with analogous definitions for left weak Hopf monads and weak Hopf
monads. 

Suppose that $T$ is a weak bimonad on a Cauchy complete monoidal category
$\M$, and $R$ the corresponding separable Frobenius monoid, with forgetful
functor $V:{}_R\M_R\to\M$ and $G\dashv V$. 
To say that $\widetilde{T}$ is right Hopf is to say that for all 
$\widetilde{X},\widetilde{Y}\in{}_R\M_R$, the canonical morphism 
$$\xymatrix{\relax
\widetilde{T}(\widetilde{T}\widetilde{X}\ox_R \widetilde{Y})
\ar[r]^{\widetilde{\can}_{\widetilde{X},\widetilde{Y}}} & \relax  
\widetilde{T}\widetilde{X}\ox_R\widetilde{T}\widetilde{Y} }$$
is invertible. Since every $\widetilde{X}\in{}_R\M_R$ is (naturally) a retract
of  one of the form $GX$, this will be the case precisely when 
$$\xymatrix{\relax
\widetilde{T}(\widetilde{T}GX\ox_R GY) \ar[r]^-{\widetilde{\can}_{GX,GY}} &
\relax \widetilde{T}GX\ox_R\widetilde{T}GY}$$ 
is invertible. Now
$$\widetilde{T}(\widetilde{T}GX\ox_R GY) \cong \widetilde{T}\bigl(TX\ox_R
(R\ox Y\ox R)\bigr) \cong \widetilde{T}(TX\ox Y\ox R)$$
and
 $V\widetilde{T}(TX\ox Y\ox R)$ is a retract of $T(TX\ox Y)$
by construction of $\widetilde{T}$, while 
$$V(\widetilde{T}GX\ox_R\widetilde{T}GY)\cong  TX\Box TY $$
which is a retract of $TX\ox TY$. Thus we obtain a composite map
\begin{equation}\label{eq:tilde-can}
\xymatrix{
T(TX\ox Y) \ar@{>>}[r]^-{q_{X,Y}} & \relax
V\widetilde{T}(\widetilde{T}GX\ox_R GY) \ar[r]^{\widetilde{\can}_{GX,GY}} & \relax
V(\widetilde{T}GX\ox_R\widetilde{T}GY)~ \ar@{>->}[r]^-{i_{TX,TY}}  & TX\ox TY }
\end{equation}
which turns out to be the canonical map $\can_{X,Y}$ associated to $T$ itself.

Now the inclusion $TX\Box TY\to TX\ox TY$ is the section for a splitting
of the idempotent $E_{TX,TY}$ on $TX\ox TY$ defined in \eqref{eq:E_AB}.

On the other hand, the quotient $T(TX\ox Y)\to
V\widetilde{T}(\widetilde{T}GX\ox_R GY)$ 
is the retraction of a splitting of an idempotent $F_{X,Y}$ on $T(TX\ox Y)$ 
defined by 
\begin{equation}
\label{eq:F_XY}
\xymatrix@C=2pc{
T(TX \mox Y) \ar[rr]^-{T(\delta\circ\eta \mox TX \mox Y \mox \eta)}&&
T(R \mox R\mox TX \mox Y \mox R)
\ar[rr]^-{T(R \mox \beta_{TX} \mox Y \mox R)}&&
T(R \mox TX \mox Y \mox R)\ar[r]^-{\lambda_{TX \mox Y}}&
T(TX \mox Y), 
}
\end{equation}
where $\beta_{TX}$ denotes the left $R$-action on $TX$ and $\lambda$ is the
natural transformation in \eqref{eq:q_X}. 

To say that $\widetilde{\can}$ is invertible, is to say that $\can$ induces 
an isomorphism between the splittings of the idempotents $F_{X,Y}$ and 
$E_{TX,TY}$. We then call $T$ a weak right Hopf monad:

\begin{definition}
\label{def:whm}
A weak bimonad $T$ on a monoidal category $(\M,\ox, K)$ is said to be a {\em
  weak right Hopf monad} provided that there is a natural transformation
$\chi_{X,Y}:TX \ox TY \to T(TX \ox Y)$ 
such that, for the canonical natural transformation $\can$ of $T$ in
\eqref{eq:can}, for the idempotent morphisms $E_{TX,TY}$ and
$F_{X,Y}$ \eqref{eq:F_XY}, and for any objects $X,Y$ of $\M$, 
\begin{equation}
\label{eq:whm}
\chi_{X,Y} \circ E_{TX,TY} = \chi_{X,Y}= F_{X,Y} \circ \chi_{X,Y},\qquad
\chi_{X,Y} \circ \can_{X,Y}= F_{X,Y},\qquad
\can_{X,Y}\circ \chi_{X,Y} = E_{TX,TY}.
\end{equation}
\end{definition}

The definition just given makes sense for any monoidal category \M, 
but is motivated by the following theorem, which requires \M~to be Cauchy
complete. 

\begin{theorem}
\label{thm:whm}
For any weak bimonad $T$ on a Cauchy complete monoidal category $(\M,\ox,K)$,
and the associated bimonad ${\widetilde T}$ in Remark \ref{rem:T_vs_Ttilde},
the following assertions are equivalent.
\begin{itemize}
\item[{(1)}] The canonical natural transformation $\widetilde{\can}$ 
 of ${\widetilde T}$ as in \eqref{eq:can}, is an isomorphism; that is,
 ${\widetilde T}$ is a right Hopf monad.  
\item[{(2)}] There is a natural transformation $\chi_{X,Y}:TX \ox TY
  \to T(TX \ox Y)$ obeying \eqref{eq:whm}; that is, $T$ is a
  weak right Hopf monad. 
\end{itemize}
\end{theorem}

\begin{proof}
The equations in \eqref{eq:whm} state exactly that the morphism induced by 
$\chi_{X,Y}$ between the splittings of $F_{X,Y}$ and $E_{TX,TY}$ is inverse to the
morphism $\widetilde{\can}_{X,Y}$ induced by $\can_{X,Y}$ between the splittings
of $E_{TX,TY}$ and $F_{X,Y}$.  
\end{proof}

\begin{remark}
Consider a weak right Hopf monad $T$ on a Cauchy complete monoidal category
$\M$ with corresponding separable Frobenius monoid $R$. By Theorem
\ref{thm:whm} and \cite[Theorem 3.6]{BLV} we conclude that whenever the
category of $R$-bimodules is right closed, this closed structure lifts to the
Eilenberg-Moore category $\M^T$. The category of $R$-bimodules is right closed
whenever $\M$ is right closed (in which case the internal homs are defined by splitting an appropriate idempotent natural transformation).
\end{remark}

Next we show that, as expected, a weak bimonoid $B$ in a Cauchy complete
braided monoidal category $\M$, induces a right weak Hopf monad $\bullet \ox
B$ on $\M$ if and only if it is a weak Hopf monoid in the sense of
\cite{AAatal:weak_cleft}, \cite{PaSt:WH_mon}.  

\begin{lemma}
\label{lem:can_vs_antipode}
For an arbitrary category $\C$, consider a functor $T:\C\to \C$ which admits
both a monad structure $\underline{T}=(T,m,u)$ and a comonad structure
$\overline{T}=(T,d,e)$. 
  Denote by ${\underline U}:\C^{\underline{T}} \to \C$ and by ${\overline
    U}:\C^{\overline{T}} \to \C$  the corresponding forgetful functors with
  respective left adjoint ${\underline F}:\C \to \C^{\underline{T}}$ and right
  adjoint $\overline F:\C \to \C^{\overline{T}}$. The following monoids (in
  $\mathsf{Set}$) are isomorphic.  
\begin{itemize}
\item[{(1)}] The monoid of natural transformations ${\overline
    F}{\underline U}\to {\overline F}{\underline U}$, with multiplication
  given by the composition of natural transformations.
\item[{(2)}] The monoid of those natural transformations 
$\gamma:{\overline F} T\to {\overline F} T$
for which 
${\overline F} m \circ \gamma T = \gamma \circ {\overline F} m$,
with multiplication given by the composition of natural transformations. 
\item[{(3)}] The monoid of natural transformations $T \to T$, with
  multiplication given by the `convolution product' $\varphi \ast \varphi' :=
  m\circ T\varphi'\circ \varphi T \circ d$.
\end{itemize}
\end{lemma}

\begin{proof}
(1)$\cong$(2). The stated isomorphism is given by the maps 
$\mathsf{Nat}({\overline F}{\underline U},{\overline F}{\underline U})\ni
\beta \mapsto \beta {\underline F}$, with the inverse 
$\gamma\mapsto {\overline F} {\underline U}\, \,{\underline \kappa}\circ 
\gamma {\underline U} \circ {\overline F}u{\underline U}$, where ${\underline
  \kappa}$ is the counit of the adjunction ${\underline F}
\dashv {\underline U}$.

(1)$\cong$(3). This is the adjunction isomorphism 
$\mathsf{Nat}({\overline F}{\underline U},{\overline F}{\underline U}) \cong 
\mathsf{Nat}({\overline U}\,\,{\overline F},{\underline U}\,\,{\underline F})$.
\end{proof}

For a functor $T$ as in Lemma \ref{lem:can_vs_antipode},
one may consider the so-called {\em `fusion operator'} in \cite{Str:fusion},
\begin{equation}
\label{eq:fusion}
\gamma:= \big(
\xymatrix{
T^2 \ar[r]^-{dT}&
T^3 \ar[r]^-{Tm}&
T^2
}\big).
\end{equation}
Clearly, it belongs to the monoid in Lemma \ref{lem:can_vs_antipode} (2). The
corresponding element of the isomorphic monoid in Lemma
\ref{lem:can_vs_antipode} (3) is the identity natural transformation $T\to
T$. (Hence, incidentally, Lemma \ref{lem:can_vs_antipode} provides an
alternative proof of \cite[Theorem 5.5]{MesWis:Hopf_mon}).  

\begin{lemma}
\label{lem:wh_bimonoid}
For a weak bimonoid $(B,\mu,\eta,\delta,\varepsilon)$ in a Cauchy complete
braided monoidal category $(\M,\ox,K,c)$,   and its induced weak bimonad
$T:=\bullet \ox B$, the following assertions hold, for any objects $X,Y$ of
$\M$. 
\begin{itemize}
\item[{(i)}] For the natural transformation \eqref{eq:can} of $T = \bullet \ox
  B$,
$$
\can_{X,Y}=(X \ox c_{Y,B}\ox B)\circ (X \ox Y \ox \gamma_K) \ox (X \ox
c_{Y,B}^{-1} \ox B),
$$
where $\gamma$ is the fusion operator \eqref{eq:fusion} for the monad and
comonad $\bullet \ox B$.
\item[{(ii)}] The idempotent natural transformation $E_{TX,TY}$ on $TX\ox TY$
  \eqref{eq:E_AB} satisfies 
\begin{equation}
\label{eq:E_K}
E_{TX,TY}=(X \ox c_{Y,B}\ox B)\circ (X \ox Y \ox E_{TK,TK}) \ox (X \ox
c_{Y,B}^{-1} \ox B).
\end{equation}
Moreover, $\bullet \ox E_{TK,TK}$ belongs to the monoid in Lemma
\ref{lem:can_vs_antipode} (2) and the corresponding element of the isomorphic
monoid in Lemma \ref{lem:can_vs_antipode} (3) is $\bullet \ox t$, where
$t$ is the composite
$$
\xymatrix@C=3pc{
B \ar[r]^-{B \ox\eta} & B^2 \ar[r]^{B\ox \delta}&
B^3 \ar[r]^-{c_{B,B}\ox B} &
B^3 \ar[r]^-{B\ox\mu} & B^2\ar[r]^-{B \ox \varepsilon}&
B.
}
$$ 
\item[{(iii)}] The idempotent natural transformation $F_{X,Y}$ on $T(TX\ox Y)$
\eqref{eq:F_XY} satisfies 
\begin{equation}
\label{eq:F_K}
F_{X,Y}=(X \ox c_{Y,B}\ox B)\circ (X \ox Y \ox F_{K,K}) \circ (X \ox
c_{Y,B}^{-1} \ox B).
\end{equation}
Moreover, $\bullet \ox F_{K,K}$ belongs to the monoid in Lemma
\ref{lem:can_vs_antipode} (2) and the corresponding element of the isomorphic
monoid in Lemma \ref{lem:can_vs_antipode} (3) is $\bullet \ox r$, where 
$r$ is the composite
$$
\xymatrix@C=3pc{
B \ar[r]^-{\eta\ox B} & B^2 \ar[r]^-{\delta\ox B}&
B^3 \ar[r]^-{B\ox c_{B,B}} &
B^3 \ar[r]^-{\mu\ox B} & B^2 \ar[r]^-{\varepsilon \ox B}&
B.
}
$$ 
\item[{(iv)}] If in addition $T:=\bullet \ox B$ is a weak right Hopf monad;
that is, there exists a natural transformation $\chi$ obeying
\eqref{eq:whm}, then   
\begin{equation}
\label{eq:chi_K}
\chi_{X,Y}=(X \ox c_{Y,B}\ox B)\circ (X \ox Y \ox \chi_{K,K}) \ox (X \ox
c_{Y,B}^{-1} \ox B)
\end{equation}
and $\bullet\ox \chi_{K,K}$ belongs to the monoid in Lemma
\ref{lem:can_vs_antipode} (2). 
\end{itemize}
\end{lemma}

\begin{proof}
Assertion (i) is immediate by relation \eqref{eq:w_bimonoid_tau} between the
opmonoidal structure $\tau_{X,Y}$ of $T= \bullet \ox B$ and the
comultiplication $\delta=\tau_{K,K}$ in $B=TK$.

(ii). Equation~\eqref{eq:E_K} follows from the formula \eqref{eq:E_AB} for
$E_{TX,TY}$ and $E_{TK,TK}$.
Then the morphism 
$$
E_{TK,TK}=\big(
\xymatrix@C=4pc{
B^2 \ar[r]^-{B^2 \ox \eta}& B^3 \ar[r]^-{B^2\ox\delta} & 
B^4 \ar[r]^-{B \ox c_{B,B}\ox B}&
B^4 \ar[r]^-{\mu \ox \mu}&
B^2
}\big)
$$
renders commutative the first diagram in  
\begin{equation}\label{eq:E_properties}
\xymatrix@C=4pc{
B^3 \ar[r]^-{B \ox E_{TK,TK}}\ar[d]_-{\mu \ox B}&
B^3 \ar[d]^-{\mu \ox B}\\
B^2 \ar[r]_-{E_{TK,TK}}&
B^2
}
\qquad
\xymatrix@C=4pc {
B^2 \ar[d]_-{B \ox \delta}\ar[r]^-{E_{TK,TK}}&
B^2 \ar[d]^-{B \ox \delta}\\
B^3 \ar[r]_-{E_{TK,TK}\ox B}&
B^3
}
\end{equation} 
by associativity of $\mu$.
By self-duality of the axioms of a weak bimonoid, the dual of
\eqref{eq:w_bimonoid_1} holds; that is, the first diagram in  
$$
\xymatrix@C=4pc @R=1pc{
B\ar[r]^-{\eta\ox B} \ar[d]_-{\eta \ox B}& B^2 \ar[r]^-{\delta\ox B} &
B^3 \ar[r]^-{B \ox c_{B,B}^{-1}}&
B^3 \ar[dd]^-{B\ox \mu}\\
B^2 \ar[d]_-{\delta\ox \delta}&\\
B^4 \ar[r]_-{B \ox c_{B,B}^{-1}\ox B}&
B^4 \ar[r]_-{B \ox \mu\ox B}& B^3 \ar[r]_-{B\ox\varepsilon\ox B} & 
B^2
}\qquad 
\xymatrix@C=4pc @R=3.5pc{
B^3 \ar[r]^-{E_{TK,TK}\ox B}&
B^3 \ar[d]^-{B \ox \varepsilon\ox B}\\
B^2 \ar[u]^-{B \ox \delta}\ar[r]_-{E_{TK,TK}}&
B^2
}
$$ 
commutes. Tensoring on the left with $B$ and then composing with 
$\mu\ox B$ gives commutativity of the diagram on the right.
It follows by coassociativity of $\delta$ that also the second diagram in
\eqref{eq:E_properties} commutes.
This proves that $\bullet \ox E_{TK,TK}$ belongs to the monoid in
Lemma~\ref{lem:can_vs_antipode} (2). The corresponding element of the isomorphic
monoid in Lemma \ref{lem:can_vs_antipode} (3) is $(B \ox \varepsilon) \circ
E_{TK,TK}\circ (\eta \ox B)=t$ as stated, by unitality of $\mu$. 

(iii). Similarly to part (ii), one easily checks that 
$$
F_{X,Y}= 
(X \ox c_{Y,B} \ox B) \circ
(X \ox Y \ox 
(\mu \ox \varepsilon \circ \mu \ox B)
\circ (B \ox c^{-1}_{B,B}\circ \delta \circ \eta \ox \delta)
)\circ
(X \ox c_{Y,B}^{-1}  \ox B),
$$
which proves \eqref{eq:chi_K}.
By associativity of $\mu$ and by coassociativity of $\delta$,
$$
F_{K,K}=\big(
\xymatrix@C=3pc{
B^2 \ar[r]^-{B\ox \eta\ox B}&
B^3 \ar[r]^-{B\ox \delta \ox \delta}&
B^5 \ar[r]^-{B \ox c^{-1}_{B,B}\ox B^2}&
B^5 \ar[r]^-{\mu \ox \mu \ox B}&
B^3 \ar[r]^-{B \ox \varepsilon \ox B}&
B^2
}\big)
$$
makes commute both diagrams in 
\begin{equation}
\label{eq:F&mu}
\xymatrix@C=4pc{
B^3 \ar[r]^-{B \ox F_{K,K}}\ar[d]_-{\mu \ox B}&
B^3 \ar[d]^-{\mu \ox B}\\
B^2 \ar[r]_-{F_{K,K}}&
B^2
}
\qquad
\xymatrix@C=4pc {
B^2 \ar[d]_-{B \ox \delta}\ar[r]^-{F_{K,K}}&
B^2 \ar[d]^-{B \ox \delta}\\
B^3 \ar[r]_-{F_{K,K}\ox B}&
B^3.
}
\end{equation} 
Thus $\bullet \ox F_{K,K}$ is an element of the monoid in
\ref{lem:can_vs_antipode} (2). By unitality of $\mu$ and counitality of
$\delta$, the corresponding element $(B \ox \varepsilon) \circ F_{K,K}\circ
(\eta \ox B)$ of the isomorphic monoid in Lemma \ref{lem:can_vs_antipode} (3)
is the stated morphism $r$.

(iv). In \eqref{eq:tilde-can}  we have seen  a relationship between the
canonical morphism $\can$ of $T$ and the canonical morphism
$\widetilde{\can}$ of the weakly lifted bimonad  $\widetilde{T}$.
Using this along with part (i), we deduce that
$\widetilde{\can}_{R\ox X \ox R,R\ox Y \ox R}$ is equal to 
$$
p_{TX,TY}\circ 
(X \ox c_{Y,B} \ox B)\circ 
(X \ox Y \ox 
i_{TK,TK}\circ 
\widetilde{\can}_{R\ox R,R\ox R} \circ 
q_{K,K}) \circ 
(X \ox c_{Y,B}^{-1} \ox B)\circ 
j_{X,Y}
$$
where $p_{TX,TY}$ is the epi part of the splitting of $E_{TX,TY}$, and $j_{X,Y}$ is
the mono part of the splitting of $F_{X,Y}$. 
Hence in view of \eqref{eq:E_K} and \eqref{eq:F_K},
$\widetilde{\can}^{-1}_{R\ox X \ox R,R\ox Y \ox R}$ is equal to 
$$
q_{X,Y}\circ 
(X \ox c_{Y,B} \ox B)\circ 
(X \ox Y \ox 
j_{K,K} \circ 
\widetilde{\can}^{-1}_{R\ox R,R\ox R} \circ 
p_{TK,TK}) \circ 
(X \ox c_{Y,B}^{-1} \ox B)\circ 
i_{TX,TY}.
$$
Thus for $\chi_{X,Y}= j_{X,Y}\circ \widetilde{\can}^{-1}_{R\ox X \ox
  R,R\ox Y \ox R}  \circ p_{TX,TY}$, the required condition \eqref{eq:chi_K}
holds. 

We need to show that $\chi_{K,K}$
induces a natural transformation as in Lemma \ref{lem:can_vs_antipode} (2). 
By part (i), $\can_{K,K}=\gamma_K$ induces such a natural
transformation. Hence in view of \eqref{eq:tilde-can}, 
since $i_{TX,TY}$ is a morphism of left $B$-modules and of right
$B$-comodules, and by \eqref{eq:F&mu},
\begin{eqnarray*}
&&(\mu\ox B)\circ 
(B \ox j_{K,K})\circ 
(B\ox \widetilde{\can}^{-1}_{R\ox R,R\ox R})=
j_{K,K}\circ \widetilde{\can}^{-1}_{R\ox R,R\ox R} \circ 
(\mu \ox_R B) \qquad \textrm{and}\\
&&(B \ox \delta) \circ 
j_{K,K}\circ 
\widetilde{\can}^{-1}_{R\ox R,R\ox R} =
(j_{K,K}\ox B)\circ 
(\widetilde{\can}^{-1}_{R\ox R,R\ox R}\ox B)\circ 
(B \ox_R \delta).
\end{eqnarray*}
Since $p_{TK,TK}$ is a
morphism of left $B$-modules and of right $B$-comodules, this implies that
$\chi_{K,K}= j_{K,K}\circ \widetilde{\can}^{-1}_{R\ox R,R\ox R} \circ
p_{TK,TK}$ belongs to the monoid in Lemma \ref{lem:can_vs_antipode} (2). 
\end{proof}

\begin{theorem}\label{thm:weakHopf}
For a weak bimonoid $B$ in a Cauchy complete braided monoidal category
$(\M,\ox,K,c)$, 
the induced functor $\bullet \ox B$ is a weak right Hopf monad if and only if
$B$ is a weak Hopf monoid. 
\end{theorem}

\begin{proof}
By Lemma \ref{lem:wh_bimonoid}, $\bullet \ox B$ is a weak right Hopf monad if
and only if (using the same notation in the lemma) there is an element 
$X \ox Y \ox \chi_{K,K}: X \ox Y \ox B \ox B \to X \ox Y \ox B \ox
B$ of the monoid in Lemma \ref{lem:can_vs_antipode} (2), such that 
$$
\chi_{K,K}\circ E_{TK,TK} = \chi_{K,K} = F_{K,K}\circ \chi_{K,K},\quad
\chi_{K,K}\circ \can_{K,K} = F_{K,K}, \quad
\can_{K,K}\circ\chi_{K,K}= E_{TK,TK}.
$$
By Lemma \ref{lem:can_vs_antipode}, this is equivalent to the existence of a
morphism $\nu:B \to B$, such that
$$
\nu \ast r = \nu = t \ast \nu, \qquad \qquad \qquad
\nu \ast B = t,\qquad \qquad \qquad
B \ast \nu = r,
$$
where the morphisms $t,r: B \to B$ are introduced in Lemma
\ref{lem:wh_bimonoid} and $\ast$ denotes the convolution product $f \ast g=\mu
\circ (f \ox g)\circ \delta$, for any morphisms $f,g:B \to B$ in $\M$.
\end{proof}

Finally we turn to connections between right weak Hopf monads and left
weak Hopf monads.
Conditions \eqref{eq:tau_1}-\eqref{eq:tau_5}
are invariant under replacing the monoidal product $\ox$ with the opposite
product ${\overline \ox}$. That is, if $(T,m,u,\tau_0, \tau)$ is a weak bimonad
on a monoidal category $(\M,\ox,K)$, then $(T,m,u,{\overline \tau}_0,
{\overline \tau})$ is a weak bimonad on $(\M,{\overline \ox},K)$, where
${\overline \tau}_0= \tau_0:TK \to K$ and ${\overline \tau}_{X,Y}= 
\tau_{Y,X}:T(X {\overline \ox} Y)=T(Y \ox X) \to TY \ox TX = TX {\overline
  \ox} TY$. 
We say that a weak bimonad $(T,m,u,\tau_0, \tau)$ is a {\em left weak Hopf
  monad} on a monoidal category $(\M,\ox,K)$ provided that $(T,m,u,{\overline
  \tau}_0, {\overline \tau})$ is a right weak Hopf monad on $(\M,{\overline
  \ox},K)$. 
Clearly, this means that the {\em left canonical map}
$$
\xymatrix{
T(X\ox TY)\ar[r]^-{\tau_{X,TY}}&
TX \ox T^2Y \ar[r]^-{TX \ox m_Y}&
TX \ox TY
}
$$
induces an isomorphism between the retracts of $T(X\ox TY)$ and $TX\ox TY$
defined as above.

Some known facts about right weak Hopf monads immediately translate to 
left weak Hopf monads:
obviously, for a weak Hopf monoid $(B,\mu,\eta,\delta,\varepsilon,\nu)$ in a 
braided monoidal category $(\M,\ox,K,c)$, the same data 
$(B,\mu,\eta,\delta,\varepsilon,\nu)$ describe a weak Hopf monoid in 
$(\M,{\overline  \ox},K,{\overline  c})$, where the braiding is given by 
${\overline  c}_{X,Y}=c_{Y,X}:X {\overline \ox }Y = Y \ox X \to X \ox Y = Y
{\overline \ox} X$.
From Theorem~\ref{thm:weakHopf} we deduce
\begin{proposition}
For a weak bimonoid $B$ in a Cauchy complete braided monoidal category
$(\M,\ox,K,c)$, the following assertions are equivalent.
\begin{enumerate}
\item the weak bimonad $\bullet\ox B$ on $(\M,\ox,K)$ is a right weak Hopf
  monad; 
\item the weak bimonad $\bullet\ox B= B {\overline  \ox} \bullet$ on 
$(\M,{\overline\ox},K)$ is a left weak Hopf  monad;  
\item $B$ is a weak Hopf monoid in $(\M,\ox,K,c)$;
\item $B$ is a weak Hopf monoid in $(\M,{\overline  \ox},K,{\overline  c})$. 
\end{enumerate}
\end{proposition}
In particular, the equivalence of the second and fourth assertions says that a
weak bimonoid in a Cauchy complete braided monoidal category is a weak
Hopf monoid if and only if it induces a left weak Hopf monad by tensoring on
the left. 

Our next aim is to describe those weak bimonoids
$(B,\mu,\eta,\delta,\varepsilon)$ in a Cauchy complete braided monoidal
category $(\M,\ox,K,c)$ for which $\bullet\ox B$ is both a right and a left
weak Hopf monad.  

Consider the weak bimonoid $B^{op}$ in $(\M,\ox,K,c^{-1})$
with the same comonoid structure $(\delta,\varepsilon)$ of $B$, multiplication
$\mu^{op}:=\mu \circ c^{-1}_{B,B}$ and unit $\eta$. Observe that via
$c_{\bullet,B}:\bullet\ox B \to B \ox \bullet$, the weak bimonads
$\bullet \ox B$ and $B^{op} \ox \bullet$ are isomorphic.
Hence the following assertions on $B$ are equivalent: 
\begin{enumerate}
\item the weak bimonad $\bullet\ox B$  on $(\M,\ox,K)$ is a left weak Hopf
  monad; 
\item the weak bimonad $B^{op} \ox \bullet$ on $(\M,\ox,K)$ is a left weak
  Hopf monad; 
\item $B^{op}$ is a weak Hopf monoid in $(\M,\ox,K,c^{-1})$;
\item there is a morphism $\nu^{op}:B \to B$ (the {\em antipode} for $B^{op}$)
such that the following diagrams commute. 
$$
\xymatrix@=1pc{
B \ar[r]^-{\delta} \ar[d]_-{\eta \ox B}& 
B^2 \ar[r]^-{\nu^{op} \ox B}&
B^2 \ar[d]^-{c^{-1}_{B,B}}\\
B^2 \ar[d]_-{\delta \ox B}&
&
B^2 \ar[d]^-{\mu}\\
B^3 \ar[r]_-{B \ox \mu}&
B^2 \ar[r]_-{B \ox \varepsilon}&
B
}\qquad
\xymatrix@=1pc{
B \ar[r]^-{\delta} \ar[d]_-{B \ox \eta}& 
B^2 \ar[r]^-{B \ox \nu^{op}}&
B^2 \ar[d]^-{c^{-1}_{B,B}}\\
B^2 \ar[d]_-{B\ox \delta}&
&
B^2 \ar[d]^-{\mu}\\
B^3 \ar[r]_-{\mu\ox B}&
B^2 \ar[r]_-{\varepsilon\ox B}&
B
}
\qquad
\xymatrix@=1pc{
B \ar[rr]^-{\delta^2} \ar[dd]_-{\nu^{op}}&&
B^3 \ar[d]^-{\nu^{op}\ox B \ox \nu^{op}}\\
&&
B^3 \ar[d]^-{(\mu\circ c_{B,B}^{-1})^2}\\
B \ar@{=}[rr]&&
B 
}
$$
\end{enumerate}
We shall use the notations $s$, $r$, $t$ in
weak Hopf monoids, as in \cite{PaSt:WH_mon}; the forms of $t$ and $r$
are recalled in Lemma~\ref{lem:wh_bimonoid} above.
The left-bottom path in the first diagram in (4) above (playing the role of
$t^{op}$) is equal to $s$. The left-bottom path in the second diagram is
conveniently denoted by $r^{op}$. The four morphisms $s$, $t$, $r$, $r^{op}$
obey the following four equations 
\begin{equation}\label{eq:appendix}
\nu\circ s =r\qquad
\nu \circ r^{op} = t\qquad
s\circ \nu =t\qquad
r^{op}\circ \nu =r.
\end{equation}
The first one is (15) in Appendix B of \cite{PaSt:WH_mon},
and the others are proved by similar steps.

Finally we are ready to provide the desired characterization:

\begin{theorem}
For a weak bimonoid $B$ in a Cauchy complete braided monoidal category
$(\M,\ox,K,c)$, the following conditions are equivalent:
\begin{enumerate}
\item the weak bimonad $\bullet\ox B$ on $(\M,\ox,K)$ is both a right and a
  left weak Hopf monad; 
\item $B$ is a weak Hopf monoid in $(\M,\ox,K,c)$ and $B^{op}$ is a weak Hopf
  monoid in $(\M,\ox,K,c^{-1})$; 
\item $B$ is a weak Hopf monoid in $(\M,\ox,K,c)$  with an invertible
  antipode $\nu$;
\item $B^{op}$ is a weak Hopf monoid in $(\M,\ox,K,c^{-1})$ with an invertible
  antipode $\nu^{op}$.
\end{enumerate}
In case (3), $\nu^{-1}$ will be an antipode for $B^{op}$; in case (4), 
$(\nu^{op})^{-1}$ will be an antipode for $B$.
\end{theorem}

\begin{proof}
We have already seen that (1) and (2) are equivalent. We show that (3) is 
equivalent to (2); the equivalence of (4) and (2) is similar.

Assume first that the (3) holds: $B$ is a weak Hopf
monoid in $(\M,\ox,K,c)$  with an invertible antipode $\nu$. In order to see
that $\nu^{-1}$ provides an antipode for the weak Hopf monoid $B^{op}$,
compose with $\nu^{-1}$ on the left the antipode axioms for $B$. The first two 
antipode axioms for $B^{op}$ follow from the respective axiom for $B$, using
the anti-multiplicativity of $\nu$ \cite[(17)]{PaSt:WH_mon} and the
identities $\nu \circ s=r$ and $\nu \circ r^{op}=t$, respectively. 
The third antipode axiom for $B^{op}$ follows from the
corresponding axiom for $B$ by \cite[(17),(6)]{PaSt:WH_mon}. Thus (2) holds.

Conversely, assume that (2)  holds: $B$ admits an
antipode $\nu$ and $B^{op}$ admits an antipode $\nu^{op}$. In order to see
that $\nu^{op}$ is a left inverse of $\nu$, use associativity of the
multiplication, and coassociativity of the comultiplication in $B$ to compute
the convolution product 
$
(\nu^{op}\circ \nu)\ast \nu \ast B =
\mu^2 \circ \big((\nu^{op}\circ \nu)\ox \nu \ox B \big) \circ \delta^2
$
in two different ways. On one hand, 
$$
\big((\nu^{op}\circ \nu)\ast \nu \big) \ast B 
=(\mu \circ c_{B,B}^{-1}\circ (B\ox \nu^{op})\circ \delta \circ \nu) \ast B 
=(r^{op}\circ \nu) \ast B 
= r\ast B 
= B.
$$
The first equality follows by anti-comultiplicativity of $\nu$,
cf. \cite[(16)]{PaSt:WH_mon}. The second equality is a consequence of one of
the antipode axioms for $B^{op}$. The third equality follows by the identity
$r^{op}\circ \nu = r$ \eqref{eq:appendix} and the last equality is easily derived
from the form of $r$ and axiom (b) in \cite{PaSt:WH_mon}.
On the other hand, 
$$
(\nu^{op}\circ \nu)\ast \big( \nu \ast B \big)
=(\nu^{op}\circ \nu)\ast t
= (\nu^{op}\circ \nu)\ast (s \circ \nu)
= \mu \circ c_{B,B}^{-1}\circ (s \ox \nu^{op})\circ \delta \circ \nu 
= \nu^{op} \circ \nu.
$$
The first equality follows by one of the antipode axioms for $B$. The second
equality follows by the identity $s\circ \nu =t$ \eqref{eq:appendix} and the
third one  follows by anti-comultiplicativity of $\nu$,
cf. \cite[(16)]{PaSt:WH_mon}. The last equality follows by the weak Hopf
monoid identity $t\ast \nu =\nu$ applied to $B^{op}$.
A symmetrical reasoning shows that $\nu^{op}$ is also a right inverse of
$\nu$: 
By \cite[(17)]{PaSt:WH_mon}, one of the antipode axioms for $B^{op}$, the
identities $\nu \circ s =r$ \eqref{eq:appendix} and $r\ast B=B$, 
$$
\big((\nu\circ \nu^{op})\ast \nu \big) \ast B =B.
$$
On the other hand, by one of the antipode axioms for $B$, the identity $\nu
\circ r^{op}=t$ \eqref{eq:appendix}, by \cite[(17)]{PaSt:WH_mon} and the weak
Hopf monoid identity $\nu \ast r =\nu$ applied to $B^{op}$,
$$
(\nu\circ \nu^{op})\ast \big( \nu \ast B  \big) =\nu\circ \nu^{op}.
$$
Thus (3) holds.
\end{proof}

\section*{Acknowledgements}
\noindent
The authors gratefully acknowledge partial support from the Australian
Research Council (project DP0771252). GB thanks also for a support of the
Hungarian Scientific Research Fund (OTKA K68195). She is grateful to the members
of the Mathematics Department and of the Department of Computing at Macquarie
University for an inspiring atmosphere and a warm hospitality experienced
during her visit in Jan-Feb 2010.  


\end{document}